\documentclass[%
 aip,
 cha,
 floatfix,
 amsmath,amssymb, 
 reprint,%
]{revtex4-1}

\usepackage{graphicx}
\usepackage{dcolumn}
\usepackage{bm}

\usepackage[utf8]{inputenc}
\usepackage[T1]{fontenc}
\usepackage{mathptmx}
\usepackage{etoolbox}

\usepackage{amsfonts, bbold, dsfont, graphicx, cancel, amsthm, array, mathtools}
\usepackage{tikz, listofitems, xcolor, bigints}

\usepackage{textcomp}
\usepackage{xcolor}
\usepackage{tabularx}

\usepackage{nicefrac}
\def\x{\bm x}
\def\z{\bm z}
\def\f{\bm f}
\def\y{\bm y}
\def\z{\bm z}
\def\w{\bm w}
\def\ftilde{\Tilde{\bm {f}}}

\def\X{\bm X}
\def\Y{\bm Y}
\def\R{\mathds{R}}

\definecolor{aquamarine}{rgb}{0.5, 1.0, 0.83}
\colorlet{myred}{red!80!black}
\colorlet{myblue}{blue!80!black}

\usetikzlibrary{automata, positioning}

\tikzstyle{data}=[thick, rectangle, draw=aquamarine, fill=aquamarine!10, minimum size=45,inner sep=0.5,outer sep=0.6]
\tikzstyle{neuron}=[thick, circle, draw=myblue!30, fill=myblue!8, minimum size=50, inner sep=0.5, outer sep=0.6]
\tikzstyle{layer}=[thick, rounded corners, draw=myred!30, fill=myred!10, minimum size=40, inner sep=0.5, outer sep=0.6]
\tikzstyle{connect}=[thick, gray]
\tikzstyle{connect arrow}=[->,thick,gray,shorten <=0.5,shorten >=1]
\tikzstyle{annot} = [text width=5em, text centered]
\tikzstyle{annot2} = [text width=7em, text centered]

\makeatletter
\def\@email#1#2{%
 \endgroup
 \patchcmd{\titleblock@produce}
  {\frontmatter@RRAPformat}
  {\frontmatter@RRAPformat{\produce@RRAP{*#1\href{mailto:#2}{#2}}}\frontmatter@RRAPformat}
  {}{}
}%
\makeatother
\begin{document}

\preprint{AIP/123-QED}

\title[Symbolic Regression via Neural Networks]{Symbolic Regression via Neural Networks}
\author{N. Boddupalli}
 \altaffiliation[Equal Contribution]{}
\author{T. Matchen}%
 \altaffiliation[Equal Contribution]{}
 \email{nibodh@ucsb.edu.}
\author{J. Moehlis}
\affiliation{%
$^{1)}$Department of Mechanical Engineering, University of California, Santa Barbara, CA 93106, USA.
}%

\date{Chaos 33, 083150 (2023)}

\begin{abstract}
Identifying governing equations for a dynamical system is a topic of critical interest across an array of disciplines, from mathematics to engineering to biology. Machine learning -- specifically deep learning -- techniques have shown their capabilities in approximating dynamics from data, but a shortcoming of traditional deep learning is that there is little insight into the underlying mapping beyond its numerical output for a given input. This limits their utility in analysis beyond simple prediction. Simultaneously, a number of strategies exist which identify models based on a fixed dictionary of basis functions, but most either require some intuition or insight about the system, or are susceptible to overfitting or a lack of parsimony. Here we present a novel approach that combines the flexibility and accuracy of deep learning approaches with the utility of symbolic solutions: a deep neural network that generates a symbolic expression for the governing equations. We first describe the architecture for our model, then show the accuracy of our algorithm across a range of classical dynamical systems.
\end{abstract}

\maketitle

\begin{quotation}
The dynamics of quantities of interest are widely modeled as differential equations, often derived from first principles.  However, this is not always possible, especially when the underlying mechanisms are unknown or complex. The identification of models from data has seen significant advances with the advent of machine learning. While deep neural networks have enabled sufficient accuracy in forecasting dynamic data with unprecedented versatility, the models they represent lack closed-form expressions that can be conducive to interpretation and analysis. Here we present an algorithm that identifies parsimonious closed-form ordinary differential equations from noisy data using a novel deep learning architecture and an  information criterion.
\end{quotation}

\section{INTRODUCTION}

Mathematical models for systems of scientific, technological, and societal interest have been pursued for centuries.  Models which capture a system's dynamics are of particular importance because they allow prediction and analysis of how quantities of interest change with time and can often be used to develop control algorithms.  For many systems, it is possible to derive mathematical models from first principles such as Newton's laws, Maxwell's equations, or the Navier-Stokes equations.  Such models have had far-reaching success over a broad range of length and time-scales, and have led to incredible advances in fields ranging from fluid dynamics to solid mechanics to robotics and beyond.
Unfortunately, deriving models from first principles is not always feasible,
which has led researchers to explore ways to deduce mathematical models from observed data, a process known as system identification~\cite{ljun99}.   
Here we briefly summarize several notable methods for system identification~\cite{brunton2019data}.

A number of authors have approached system identification by fitting coefficients of a linear combination of basis functions, dating at least back to Crutchfield and McNamara~\cite{crut87}.  The set of basis functions typically includes nonlinear terms, for example terms which would arise in a Taylor series expansion about the origin of the system~\cite{crut87, yao07, wang11, lai21} or a broader class of functions~\cite{brun16}.  The coefficients of the basis functions are determined through comparison of the original data points with points from computed solutions to the fitted models.  Various sparse optimization approaches have been proposed to try to avoid an unwieldy number of terms in the resulting models~\cite{napo08, wang11, brun16,lai21}.  Given its recent popularity, here we highlight the Sparse Identification of Nonlinear Dynamics (SINDy) algorithm~\cite{brun16}, in which the right-hand side of the differential equation model is found as the optimal sparse linear combination of user-specified candidate functions, which may be nonlinear.  The SINDy algorithm has been successfully applied to a variety of models 
including the Lorenz equations and vortex shedding past a cylinder
~\cite{brun16}, and has been extended to allow for rational function nonlinearities which commonly appear in biological networks~\cite{mang16}. It was also combined with deep learning to determine which terms are needed to accurately capture the dynamics seen in  data~\cite{champion2019data}. A limitation of all such approaches is that they require pre-specified basis functions whose linear combination gives the model.  We call such methods ``fixed dictionary'' methods.  Mathematically, fixed dictionary system identification approaches find a model
\begin{equation}\label{fixed_dict}
    \frac{dx}{dt} = \sum_i a_i \phi_i(x)
    \end{equation}
by making smart choices for the $a_i$'s, where $\{\phi_i(x)\}$ is a pre-specified fixed dictionary of functions.

Also of recent interest is dynamic mode decomposition (DMD), which discovers coherent spatiotemporal modes from measurements of complex systems, and has connections to Koopman operator theory~\cite{mezic2005spectral, budi12, arba17}.  DMD was initially developed in the fluid dynamics community~\cite{schm10}, and has since been applied to many other systems 
~\cite{kutz16}.  A type of DMD known as extended dynamic mode decomposition (EDMD) is a regression onto locally linear dynamics of a fixed dictionary of candidate functions, and can be used to generate a nonlinear closed-form mathematical model for a system based on data~\cite{edmd}.  However, EDMD can face ``curse of dimensionality'' challenges and typically requires a large amount of data to be accurate.

In another class of system identification methods, which we call ``generative dictionary'' methods, one only pre-specifies a set of primitive operations and functions, and the method creates the dictionary of possible model terms from combinations and compositions of these primitives.  Historically, such methods have been referred to as ``symbolic regression'', because they search a space of mathematical expressions to find a model, as opposed to conventional regression techniques which optimize parameters for a pre-specified model structure (i.e., a fixed dictionary).  Typical symbolic regression approaches to system identification~\cite{bong07,quad16,quade2019glyph} learn models from observed data by randomly combining various terms and operations, and using genetic programming to ``mutate'' the candidate solutions according to a fitness-weighted selection mechanism. One such approach~\cite{schmidt2009distilling} was to restrict the space of possible models by searching for invariant sets. While these invariants involved estimating higher order time and space derivatives that grow combinatorially in number with increase in dimensions, it was also only suitable to model systems that obeyed conservation laws. Unfortunately, genetic programming strategies typically suffer from bloat~\cite{banz02}, excessive terms, and overly complex representations of systems, making identification of the most physically meaningful representation of the data difficult.  

We emphasize that there is an important distinction between fixed dictionary and generative dictionary system identification methods: for the former, only terms which have been included in the original fixed dictionary can appear in the model.  For the latter, any term which can be generated from combinations and compositions of the primitives can appear.  In particular, let $S_k(x)$ be a (possibly $k$-dependent) subset of all possible functions of $x$ generated by the primitive operations and functions.  Mathematically, a generative dictionary system identification method seeks to find a model
\begin{equation}\label{generative_dict}
\frac{dx}{dt} = F(S_K \circ \cdots \circ S_2 \circ S_1),
\end{equation}
by making smart choices for which function $F$ and subsets $S_k$ to choose.  Here the composition $\circ$ is interpreted in an element-wise manner; for example, $S_2 \circ S_1$ can contain functions $f_2(f_1(x))$ for any $f_1 \in S_1$ and $f_2 \in S_2$.  Moreover, each $S_k$ includes the Identity operator, so, for example, $S_2 \circ S_1$ can contain any function in $S_1$. The generative dictionary is given by the set $S_K \circ \cdots \circ S_2 \circ S_1$, and both the generative dictionary and the function $F$ are determined as part of the system identification algorithm.

An alternative approach to system identification is to use deep neural networks (DNNs) to obtain a ``black box'' model in which the current state of a system is mapped to the future state without knowledge of the inner workings of the transformation~\cite{lape88,kube92} \cite{timeseriesreview}.
For system identification, it is notable that DNNs do not make assumptions about the terms in a model, unlike SINDy or EDMD.  This comes at the cost of estimating the network weights iteratively, but fast optimization algorithms benefit from today's computational power that is higher and cheaper than before. On the other hand, 
generalizability of the obtained weights to unseen data-sets even in the same regime is not guaranteed. 

In this paper, we present a novel neural network framework for identifying interpretable models from time series or vector-field data without {\it a priori} knowledge of the governing dynamics. The models are built up from primitive operations like addition, multiplication, exponentiation, and other user-defined functions of state variables that commonly occur in dynamical systems. Drawing inspiration from DNNs, we allow the breadth of the network to determine the number of allowable coefficients in a certain set of functions, and the depth to determine compositions among them. Our approach allows the generation of interpretable functions of state variables including linear combinations, polynomials, and non-polynomial functions such as fractional and negative powers, logistic functions, and expressions that might arise in chemical kinetics or a conductance-based model for neural activity. The complexity of these generated functions that appear in the model is determined directly by the depth and width of the neural network rather than having to pre-specify them exactly. Our neural network architecture is novel in how the weights of the network determine the form of the terms in the generated dictionary, which allows a wider variety of model terms than other recent system identification approaches that use neural networks~\cite{kim21, fronk23}.  Moreover, the number of terms in the generated dictionary remains small enough to handle computationally, but large enough to capture a rich set of possibilities. Available algorithms for optimizing DNNs are used to optimize the coefficients of state-variables in these expressions. To obtain parsimonious models, we use a combination of sparsity-inducing regularization and the Akaike Information Criterion (AIC)~\cite{akaike1974new}, which selects between candidate models, obtained through perturbations of model parameters within user-specified tolerances, to balance model simplicity and accuracy. We call this SymANNTEx (pronounced as ``semantics''), for \underline{Sym}bolic, \underline{A}rtificial \underline{N}eural \underline{N}etwork-\underline{T}rained \underline{Ex}pressions.

\subsection{Preliminaries}\label{prelims}

We consider dynamical systems of the form
\begin{equation}\label{sys}
    \dot{\x} = \f(\x),
\end{equation}
where $\x \in \R^n$ and $\f: \R^n \to \R^n$. The $i^{\rm th}$ state of $\x$ will be denoted as $x_i$. We consider data-points sampled discretely -- but not necessarily from a single trajectory -- in time, and denote the $j^{\rm th}$ data-point as $\x^{[j]}$. Given $m$ such data-points $\{\x^{[j]}\}_{j = 1}^m$, we train our neural network on data that approximates $\f(\x^{[j]}) = \dot{\x}|_{\x = \x^{[j]}} \equiv \dot{\x}^{[j]}$ at those points, as follows.

Since real data is expected to be noisy, we consider states $x_i^{[j]} + x^{rms}_i \eta_i^{\sigma_1 \ [j]}$, where $\eta_i^{\sigma_1 \ [j]}$ is independently drawn from the normal distribution $\mathcal{N}(0, \sigma_1^2)$ centered at zero and characterized by standard deviation $\sigma_1$, which is scaled by the Root Mean Square (RMS) of the corresponding state
\begin{equation*}
x^{rms}_i \equiv \sqrt{\frac{1}{m} \sum_{j = 1}^m {(x_i^{[j]})}^2},
\end{equation*}
which is the $l^2$ norm of the corresponding state averaged over the $m$ data points. We write this using the element-wise product $(\x^{rms} \odot \bm \eta^{\sigma_1})_i^{[j]} \equiv x^{rms}_i \eta_i^{\sigma_1 \ [j]}$ and denote the state data as
\begin{equation}\label{noisy_data}
\{\x^{[j]} + (\x^{rms} \odot \bm \eta^{\sigma_1})^{[j]}\}_{j = 1}^m \equiv \{\x^{[j]}_{\sigma_1}\}_{j = 1}^m \equiv \X.
\end{equation}
For example, if we consider noise drawn from a distribution of standard deviation $\sigma_1 = 0.01$, we refer to it as $1\%$ noise relative to the corresponding component of $\x^{rms}$.
Similarly, we add noise to the corresponding derivative data and define the training data $\Y$ as follows
\begin{equation}\label{training_data}
\{\dot{\x}^{[j]} + (\dot{\x}^{rms} \odot \bm \eta^{\sigma_2})^{[j]}\}_{j = 1}^m \equiv \{\dot{\x}^{[j]}_{\sigma_2}\}_{j = 1}^m \equiv \Y,
\end{equation}
where $\eta_i^{\sigma_2 \ [j]}$ are each independently drawn from $\mathcal{N}(0, \sigma_2^2)$ and $\dot{\x}^{rms}$ is the vector of RMS of derivatives at points considered which is computed as
\begin{equation*}
\dot{x}^{rms}_i \equiv \sqrt{\frac{1}{m} \sum_{j = 1}^m {(f_i(\x^{[j]}))}^2}.
\end{equation*}
We note that derivative data could alternatively be obtained from noisy state data using various numerical differentiation techniques~\cite{van22b}.  The obtained derivative data may be imprecise, which is what $\sigma_2$ is intended to mimic. 

Note that we scale the noise added to the states and the derivatives by the RMS of the respective corresponding states and derivatives because our examples span varying length and time scales, and we found that un-scaled additive noise of a certain level can have negligible effect on one example or state, while being detrimental to another.  But we have also verified that our approach successfully identifies models when (sufficiently weak) un-scaled additive noise is used.

Our network's objective is to map input (noisy state) data $\X$ to the training (noisy derivative) data $\Y$. We call $\f(\x)$ in Eq.~\ref{sys} the ``ground truth model'', the model represented by the network $\x_{\sigma_1}^{[j]} \mapsto \ftilde(\x_{\sigma_1}^{[j]})$ the ``network model'', and the model $\x_{\sigma_1}^{[j]} \mapsto \tilde{\f}_A(\x_{\sigma_1}^{[j]})$ chosen after perturbing the network model and evaluating AIC scores as the ``selected network model''. This can also be directly applied to identify maps $\x(k+1) = \bm F(\x(k))$ with a closed-form expression. Note that even the ground truth model likely gives non-zero errors because the training data $\Y$ is not exactly the derivatives evaluated at the data points $\X$ as the expected values $E(\dot{\x}_{\sigma_2}) \neq E(\f(\x_{\sigma_1}))$ for our examples even with $\sigma_1 = \sigma_2$.

\section{Network Architecture}

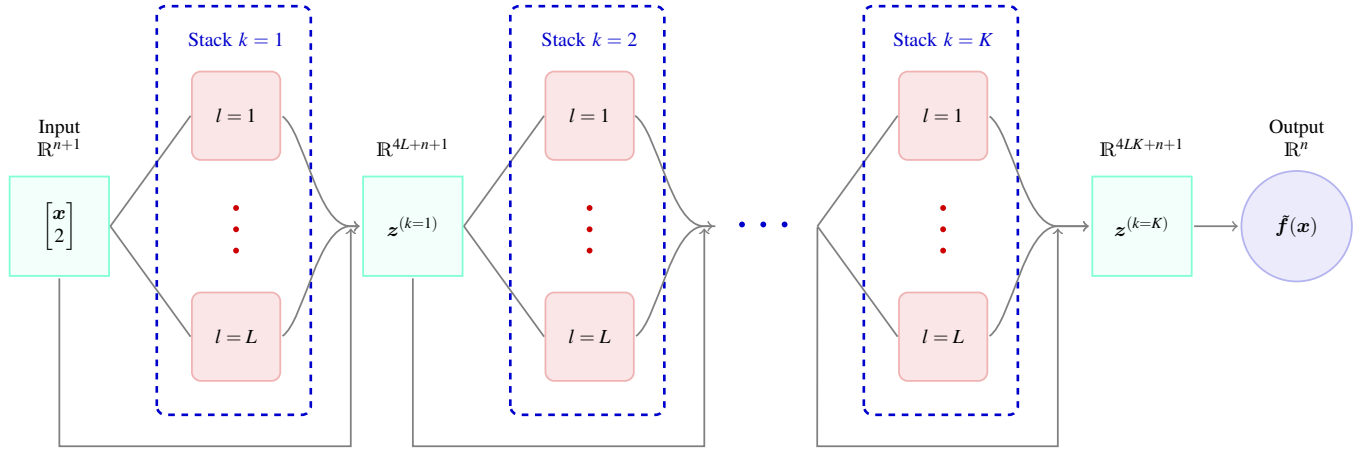
\begin{figure*}[t]
    \centering
    \resizebox{\textwidth}{!}{%
            \begin{tikzpicture}

\node[data] (input) {$\begin{bmatrix} \bm x \\ 2 \end{bmatrix}$};

\node[layer, xshift=80, yshift=-50, at=(input)] (k_1_l_1) {$l = L$};
\node[layer, xshift=80, yshift=50, at=(input)] (k_1_l_L) {$l = 1$};

\path (k_1_l_1) -- (k_1_l_L) node [myred, font=\Huge, midway, sloped] {$\dots$};

\draw[connect] (input.east) -- (k_1_l_1.west);
\draw[connect] (input.east) -- (k_1_l_L.west);

\node[data, xshift=160, yshift=0, at=(input)] (k_1) {$\z^{(k=1)}$};

\draw[rounded corners, myblue, very thick, dashed] (1.55, -3) rectangle (4, 3.5) {};

\coordinate (k_1_merge) at ([xshift=-5]k_1.west);
\coordinate (input_under) at (0, -3.5);

\coordinate (perp_point) at (input_under -| k_1_merge);
\draw[connect arrow] (input.south) -- (input_under) -- (perp_point) -- (k_1_merge);

\draw[connect arrow] (k_1_l_1.east)  .. controls (3.9, -1.5) and (4.1,0) .. (k_1_merge) -- (k_1.west);
\draw[connect arrow] (k_1_l_L.east) .. controls (3.9, 1.5) and (4.1,0) .. (k_1_merge) -- (k_1.west);

\node[annot, above = 1, align= center, at=(input)] {Input\\$\R^{n+1}$};

\node[annot2, myblue, above = 1, align= center, at=(k_1_l_L)] {Stack $k = 1$};

\node[annot, above = 1, align = center, at = (k_1)] {$\R^{4L+n+1}$};

\node[layer, xshift=80, yshift=-50, at=(k_1)] (k_2_l_1) {$l = L$};
\node[layer, xshift=80, yshift=50, at=(k_1)] (k_2_l_L) {$l = 1$};

\path (k_2_l_1) -- (k_2_l_L) node [myred, font=\Huge, midway, sloped] {$\dots$};

\draw[connect] (k_1.east) -- (k_2_l_1.west);
\draw[connect] (k_1.east) -- (k_2_l_L.west);

\node[annot2, myblue, above = 1, align= center, at=(k_2_l_L)] {Stack $k = 2$};

\draw[rounded corners, myblue, very thick, dashed, xshift= 160] (1.55, -3) rectangle (4, 3.5) {};

\coordinate (k_1_merge) at ([xshift=-5]k_1.west);

\node[data, draw = black!0, fill=black!0, xshift=320, yshift=0, at=(input)] (k_dots) {\Huge \textcolor{myblue}{$\cdots$}};

\coordinate (k_dots_merge) at ([xshift=-5]k_dots.west);

\coordinate (perp_point_dots) at (input_under -| k_dots_merge);
\draw[connect arrow] (k_1.south) -- ([xshift=160]input_under) -- (perp_point_dots) -- (k_dots_merge);

\draw[connect arrow] (k_2_l_1.east)  .. controls ([xshift=160]3.9, -1.5) and ([xshift=160]4.1,0) .. (k_dots_merge) -- (k_dots.west);
\draw[connect] (k_2_l_L.east) .. controls ([xshift=160]3.9,1.5) and ([xshift=160]4.1,0) .. (k_dots_merge) -- (k_dots.west);

\node[layer, xshift=80, yshift=-50, at=(k_dots)] (k_K_l_1) {$l = L$};
\node[layer, xshift=80, yshift=50, at=(k_dots)] (k_K_l_L) {$l = 1$};

\path (k_K_l_1) -- (k_K_l_L) node [myred, font=\Huge, midway, sloped] {$\dots$};

\draw[connect] (k_dots.east) -- (k_K_l_1.west);
\draw[connect] (k_dots.east) -- (k_K_l_L.west);

\node[annot2, myblue, above = 1, align= center, at=(k_K_l_L)] {Stack $k = K$};

\draw[rounded corners, myblue, very thick, dashed, xshift= 320] (1.55, -3) rectangle (4, 3.5) {};

\coordinate (k_dots_merge) at ([xshift=-5]k_dots.west);

\node[data, xshift=170, yshift=0, at=(k_dots)] (k_K) {$\z^{(k=K)}$};

\coordinate (k_K_merge) at ([xshift=-15]k_K.west);

\coordinate (perp_point_K) at (input_under -| k_K_merge);
\draw[connect arrow] (k_dots.east) -- ([xshift=343]input_under) -- (perp_point_K) -- (k_K_merge);

\draw[connect arrow] (k_K_l_1.east)  .. controls ([xshift=320]3.9, -1.5) and ([xshift=320]4.1,0) .. (k_K_merge) -- (k_K.west);
\draw[connect arrow] (k_K_l_L.east) .. controls ([xshift=320]3.9, 1.5) and ([xshift=320]4.1,0) .. (k_K_merge) -- (k_K.west);

\node[annot, above = 1, align = center, at = (k_K)] {$\R^{4LK+n+1}$};

\node[neuron, xshift=70, at=(k_K)] (dense_out) {$\ftilde(\x)$};

\draw[connect arrow] (k_K) -- (dense_out);

\node[annot, above = 1, align = center, at = (dense_out)] {Output\\$\R^n$};


\end{tikzpicture}
        }
        \caption{Network architecture of $K$ stacks with $L$ operational layers within each stack. Input on the far left is the state data $\x \in \X$ concatenated with the constant $2$, making the input $\R^{n+1}$. Note the  $4L(k-1) + n+ 1$ dimensional output from $(k-1)^{\text{th}}$ stack is concatenated to the $4L$ dimensional output of the $k^{\text{th}}$ stack to give a $4Lk + n+ 1$ dimensional output. After $K$ such stacks, the $4LK + n+ 1$ dimensional output of the $K^{\text{th}}$ stack is reduced to $\R^n$ to give the network model $\ftilde(\x)$ on the far right. This is optimized to fit the training (noisy derivative) data $\dot{\x}_{\sigma_2} \in \Y$ to minimize the loss shown in Eq.~\ref{loss}. The details of each operational layer (red) are further illustrated in Fig.~\ref{Layer_architecture}.}
    \label{model_architecture}
\end{figure*}

Our approach to symbolic regression fuses the adaptability of artificial neural networks with the interpretability of dictionary-based identification methods via a generative dictionary.
By design, the generated candidate functions which describe the dynamics have global support. The width of the network allows different coefficients, e.g. $\sin(a_1 x_1) + \sin(a_2 x_1)$, and the depth gives compositions of increasing complexity that have closed-form expressions, e.g. $\sin(x_1^2 x_2)$.  Specifically, the network initially takes as input the $n$ system states and an additional constant value of 2 to aid with scaling, for $n+1$ total states; notationally we let $x_{n+1} = 2$. If desired, time could be included as well, allowing for $n+2$ initial states. The neural network possesses two principal tiers of organization:  ``stacks'' and ``operational layers''. Stacks serve as the higher-level organizational structure; they modify the data in series, with the output of each stack serving as an input to subsequent stacks. The network always possesses $K\ge1$ stacks. An individual stack is in turn composed of a predetermined number of operational layers, which generate new terms for use both in subsequent stacks and the final expression. As shown in Fig.~\ref{model_architecture}, these operational layers function in parallel, sharing a common input from the original data as well as any prior stacks. Each stack possesses $L\ge1$ instances of each operational layer.

Each operational layer consists of four types of sublayers: 1) linear combinations; 2) polynomial combinations; 3) simple products of variables; and 4) common operators. As noted, if we have multiple stacks, the outputs from previous stacks serve as inputs to subsequent stacks. For example, if $K>1$, then the $k^{\text{th}}$ stack takes as input not only the initial $n+1$ states but also the $4L\left(k-1\right)$ states generated by operations in the previous stacks, which allows for highly complex analytic expressions with minimal training parameters. Each of these processes and operations is outlined below. The total number of expressions we combine to generate the final output is $4LK+n+1$ if the system is autonomous, or $4LK+n+2$ if the system is non-autonomous. These terms are combined to generate an output via a simple, fully-connected $n \times (4LK+n+1)$ dense output layer.
Weight vectors for trainable sublayers will be referred to as $\bm w$, with $w_i$ denoting the $i^{\rm th}$ component. The architecture of each operational layer is illustrated in Fig.~\ref{Layer_architecture}.

\subsection{Linear Combinations}
These sublayers are linear combinations of the pre-existing states (including $x_{n+1} = const.$)
\begin{equation}\label{linear_eq}
y = \sum_{i = 1}^{n+1} w_i x_i.
\end{equation}

\subsection{Polynomial Combinations}
The output of a polynomial combination sublayer is given by
\begin{equation}\label{poly_eq}
y = \prod_{i = 1}^{n+1} \left|x_i\right|^{w_i}.
\end{equation}
The absolute value of each state is used to ensure that each output value $y$ remains in the domain of real numbers; sign considerations are handled separately.
For example, suppose we wish to represent the equation $\dot{x_1} = {x_2}^2$ where we consider states $\left[x_1, x_2, 2\right]$; then the corresponding learned weight vector should yield $\bm w = \left[0.0, \ 2.0, \ 0.0\right]$ after training.   

\subsection{Simple Products}
In many instances, the absolute value of a state is insufficient; the simple product sublayer allows us to combine states raised only to the first power. These operations are dictated by a product rule: 

\begin{equation}\label{prod_eq}
y = \prod_{i=1}^{n+1} [\sigma\left(w_i\right)x_i + \left(1-\sigma\left(w_i\right)\right)] \equiv \prod_{i = 1}^{n+1} v_i.
\end{equation} 
Here, $\sigma\left(w_i\right) = \frac{1}{1 + e^{-w_i}}$ represents the sigmoidal activation function so that for $w_i \ll 0$, $v_i=1$ and for $w_i \gg 0$, $v_i = x_i$. 

For example, suppose $\dot{x_1} = x_1 x_2$ with states $\left[x_1, x_2, 2\right]$. Then the resulting weight vector after training could, for example, be equal to $\bm w = \left[1.5, \ 2.1,\  -0.99\right]$; here, $\sigma\left(1.5\right)\approx\sigma\left(2.1\right)\approx~1$ while $\sigma\left(-0.99\right)\approx~0$, so according to Eq.~\ref{prod_eq} yields: \begin{equation*}
\left(1 \cdot x_1+\left(1-1\right)\right)\left(1 \cdot x_2+\left(1-1\right)\right)\left(0\left(2\right)+\left(1-0\right)\right)=x_1 x_2.
\end{equation*}

\subsection{Common Operators}\label{operator_sublayer}
These sublayers are reliant on a set of common symbolic operations. Each common operator sublayer consists of two consecutive dense sublayers; the first dense sublayer returns a linear combination of the input states
\begin{equation}\label{ops_in}
y_a = \sum_{i = 1}^{n+1} w_i x_i, 
\end{equation} while the second is a linear combination of the different operations $g$ applied to the output of the first
\begin{equation} \label{ops_eq}
y = \sum_{j = 1}^{3} {g_j}\left(y_a\right) W_j. 
\end{equation} 

In the formulation of the model used to generate results in this paper, three possible operators $g$ are included: the sine function, exponential function, and sign function. The first two were selected because of their prevalence in dynamical systems across physics, chemistry, biology, and engineering, while the third is included to accommodate the absolute value requirement imposed on the polynomial sublayers. In principle, other commonly encountered functions such as saturation functions and ReLU can also be included.

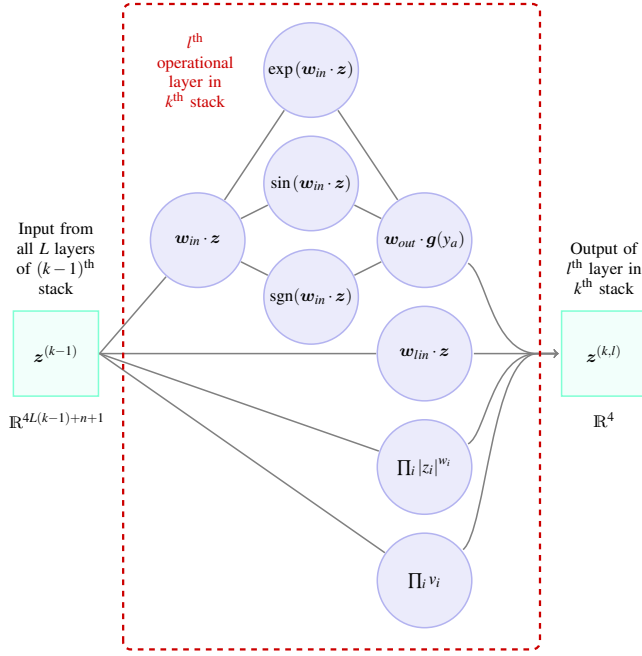
\begin{figure}[htbp]
    \centering
    \resizebox{\columnwidth}{!}{%
\begin{tikzpicture}

\node[data] (input) {$\z^{(k-1)}$};

\node[neuron, xshift=75, yshift=60, at=(input)] (ops_in) {$\w_{in} \cdot \z$};

\node[neuron, xshift=135, yshift=150, at=(input)] (exp) {$\exp{(\w_{in} \cdot \z)}$};
\node[neuron, xshift=135, yshift=90, at=(input)] (sin) {$\sin{(\w_{in} \cdot \z)}$};
\node[neuron, xshift=135, yshift=30, at=(input)] (sgn) {sgn${(\w_{in} \cdot \z)}$};

\node[neuron, xshift=195, yshift=60, at=(input)] (ops_out) {$\w_{out} \cdot \bm g(y_a)$};

\node[neuron, xshift=195, yshift=0, at=(input)] (linear) {$\w_{lin} \cdot \z$};
\node[neuron, xshift=195, yshift=-60, at=(input)] (power) {$\prod_i \left|z_i\right|^{w_i}$};
\node[neuron, xshift=195, yshift=-120, at=(input)] (product) {$\prod_i v_i$};

\node[data, xshift=290, yshift=0, at=(input)] (output) {$\z^{(k,l)}$};

\draw[connect] (input.east) -- (ops_in);

\draw[connect] (ops_in) -- (sin);
\draw[connect] (ops_in) -- (exp);
\draw[connect] (ops_in) -- (sgn);

\draw[connect] (sin) -- (ops_out);
\draw[connect] (exp) -- (ops_out);
\draw[connect] (sgn) -- (ops_out);

\draw[connect] (input.east) -- (product);
\draw[connect] (input.east) -- (power);
\draw[connect] (input.east) -- (linear);

\coordinate (output_merge) at ([xshift=-12]output.west);

\draw[connect arrow] (ops_out)  .. controls (8, 1.5) and (8.1,0) .. (output_merge) -- (output.west);
\draw[connect arrow] (linear) -- (output.west);
\draw[connect arrow] (power) .. controls (8, -1.5) and (8.1,0) .. (output_merge) -- (output.west);
\draw[connect arrow] (product) .. controls (8, -3.3) and (8.1,0) .. (output_merge) -- (output.west);

\draw[rounded corners, myred, very thick, dashed] (1.25, -5.5) rectangle (9, 6.5) {};

\node[annot, myred, left = 1.25, align=center, at=(exp)] {$l^{\text{th}}$ operational layer in $k^{\text{th}}$ stack};
\node[annot, above = 1, align= center, at=(input)] {Input from all $L$ layers of $(k-1)^{\text{th}}$ stack};
\node[annot, xshift = -2.5, yshift = -35, align = center, at = (input)] {$\R^{4L(k-1)+n+1}$};
\node[annot, above = 1, align= center, at=(output)] {Output of $l^{\text{th}}$ layer in $k^{\text{th}}$ stack};
\node[annot, xshift = 0, yshift = -35, align = center, at = (output)] {$\R^4$};

\end{tikzpicture}
        }
        \caption{Architecture of the $l^{\text{th}}$ operational layer in the $k^{\text{th}}$ stack. Input is the $4L(k-1) + n+1$ dimensional output of the $(k-1)^{\text{th}}$ stack. Output of each operational layer is $4$ dimensional, thus giving $4L$ additional states from $L$ layers in each stack. The sublayers (purple) denoted are: $\w_{lin} \cdot \z \equiv$ linear combinations (Eq.~\ref{linear_eq}), $\prod_i \left|z_i\right|^{w_i} \equiv$ polynomial combinations (Eq.~\ref{poly_eq}), $\prod_i v_i \equiv$ simple products (Eq.~\ref{prod_eq}), $\w_{in} \cdot \z, \ \w_{out}\cdot \bm g(y_a) \equiv$ linear combinations (Eqs.~\ref{ops_in}, \ref{ops_eq}) respectively where $\bm g(y_a) \equiv [\exp{(\w_{in} \cdot \z)}, \ \sin{(\w_{in} \cdot \z)}, \ \text{sgn}(\w_{in} \cdot \z)]$ in Eq.~\ref{ops_eq}. The stack and layer indices $(k,l)$ for the weights and $z$ denoted within the dashed red line are not shown in interest of clarity and the summation notation is avoided to distinguish weights sublayer-wise.
        }
    \label{Layer_architecture}
\end{figure}

We can now write down the output of each operational layer as a function. As illustrated in Fig.~\ref{Layer_architecture}, each operational layer $l$ in stack $k$ takes as input a vector in $\R^{4L(k-1)+n+1}$ and outputs a vector in $\R^4$. While we showed Eqs.~(\ref{linear_eq}-\ref{ops_eq}) for each operational layer in the first stack as functions of the input for simplicity, they can be generalized to every operational layer as
\begin{subequations}
    \begin{align}
    f^{k,l}_1(\bm z) &= \bm w_{lin}^{k, l} \cdot \bm z,\label{layer_form_1}\\
    f^{k,l}_2(\bm z) &= \prod_{i = 1}^d | z_i |^{w_{pow, i}^{k,l}},\\
    f^{k,l}_3(\bm z) &= \prod_{i=1}^d [\sigma(w_{prod}^{k,l}) z_i + (1 - \sigma(w_{prod}^{k,l}))],\\
    f^{k,l}_4(\bm z) &= \bm w_{out}^{k,l} \cdot \begin{bmatrix}\exp(\bm w_{in}^{k,l} \cdot \bm z) \\ \sin(\bm w_{in}^{k,l} \cdot \bm z) \\ \text{sgn}(\bm w_{in}^{k,l} \cdot \bm z) \end{bmatrix}.\label{layer_form_4}
\end{align}
\end{subequations}
where $d = 4L(k-1)+n+1$, $\bm z \in \R^{4L(k-1)+n+1}$. These four functions can be denoted as a vector $\bm f^{k,l}: \R^{4L(k-1)+n+1} \mapsto \R^4$ of functions:
\begin{equation}\label{layer_form}
    \bm f^{k,l}(\bm z) \equiv \begin{bmatrix}
        f^{k,1}_1 \\ f^{k,1}_2 \\ f^{k,1}_3 \\f^{k,1}_4
    \end{bmatrix}(\bm z).
\end{equation}
Now, the neural network model $\bm \ftilde(\x)$, illustrated in Fig.~\ref{model_architecture} as the output of the neural network, can be written as
\begin{equation}
    \bm \ftilde (\x) = \bm W_{out} \begin{bmatrix}
        \bm f^{K,1} \\ \bm f^{K,2} \\ \vdots \\ \bm f^{K, L} \\ \mathcal{I}
    \end{bmatrix} \circ \cdots \circ \begin{bmatrix}
        \bm f^{2,1} \\ \bm f^{2,2} \\ \vdots \\ \bm f^{2, L} \\ \mathcal{I}
    \end{bmatrix} \circ \begin{bmatrix}
        \bm f^{1,1} \\ \bm f^{1,2} \\ \vdots \\ \bm f^{1, L} \\ \mathcal{I}
    \end{bmatrix} \left( \begin{bmatrix}
        \bm x \\ 2
    \end{bmatrix} \right),
\end{equation}
where $\mathcal{I}$ is the Identity operator. Here, we can see the generative nature of our algorithm introduced in (\ref{generative_dict}), where each operational layer $l$ in stack $k$, represented as $\bm f^{k,l}$ in (\ref{layer_form}), is made up of primitives from Eqs.~(\ref{layer_form_1}-\ref{layer_form_4}) which are the sublayers elaborated earlier.

Collectively, the structure of these operational layers allows us to achieve tremendous flexibility in the expressions the network is capable of generating, while simultaneously tackling the  challenge of overfitting. While we have a significant range of possible expressions, we are able to do so with fewer trainable parameters than conventional artificial neural networks, even when the network is several stacks deep or the number of copies of the operational layers in each stack is high. Examples of the sublayer outputs when applied to a $(n+1)\times1$ vector of states are provided in Table \ref{sublayer_summary}, demonstrating how low the count of trainable parameters is when compared to fully-connected or even convolutional networks~\cite{convnets}. The smaller number of parameters is a design consequence, but offers significant advantages over conventional DNNs. 

\begin{table}[t]
\centering
\caption{Examples of sublayer outputs for $n=2$}
\begin{ruledtabular}
\begin{tabular}{ccc}
Sublayer Type & Example Output & Trainable\\
& Expression & Parameters \\
\hline
Linear Combination & $0.7x_1+1.5x_2-1$ & $n+1$ \\
Polynomial Combination & $0.25\left|x_1\right|^{-1.4}\left|x_2\right|^{0.5}$ & $n+1$\\
Simple Products & $\left(0.3+0.7x_1\right)x_2$ & $n+1$\\
Common Operators & $0.9 \sin{\left(0.2x_1-0.1 x_2\right)} \; +$ & $(n+1)+3$\\ & $0.3 \exp\left(0.2x_1-0.1x_2\right)$ & \\
\end{tabular}
\label{sublayer_summary} 
\end{ruledtabular}
\end{table}

The fact that commonly occurring terms in models are functions of state makes them analogous to activation functions in conventional DNNs. Many functions such as $\x,\ \prod_i \left|x_i\right|^{w_i}, \ \sin(\bm w \cdot \x)$, and $\ e^{\bm w \cdot \x}$ have global support, unlike conventional activation functions such as $\tanh(\bm w \cdot \x)$ and $ \ {\rm ReLU}(\bm w \cdot \x)$. While this reduces fidelity of the neural network's approximation capability, it also means that activation functions with local support need not be centered around data-points as is required in k-means clustering or other activation function centering algorithms. Consequently, far fewer parameters are required $\left(O(10^2)\right)$ compared to conventional DNNs. This carries advantages in terms of generalizability, computational economy, and data requirements. The generalizability aspect is immediate from the non-local support, meaning that newer data-points far away from the seen data $\X$ would not need activation functions centered around them. The computational requirements are far smaller both in terms of memory and processing requirements, a consequence of the smaller number of trainable parameters. The examples shown in this paper are run on a fairly standard laptop PC with an AMD Ryzen 4900HS processor and 16GB of RAM. The time taken from generating the training data to obtaining the equations and figures can range from a few minutes to an hour, depending on the dataset and neural network sizes.

\subsection{Regularization and Initialization}

We designed regularization and initialization of the neural network to achieve two primary goals: 1) Promote parsimony in the output equation and 2) Prevent exploding loss/gradient values in deeper multi-stack implementations. To this end, appropriate initializers and regularizers were selected for each layer. We will briefly detail and motivate the choices made for each operational layer here.

We make two preliminary observations before proceeding. First, parsimony is generally enforced through sparsity, but in the context of our network, zero matrix weights do not directly correspond to sparsity in most layers. There is a direct connection between parsimony and sparsity in the dense output layer and other dense layers. Second, the primary cause of exploding errors or gradients during training of multi-stack networks arises from the compounding effects of variables being multiplied by themselves, so in general initialization is designed to make most terms initially nearly constant-valued.

\subsubsection{Linear Combinations}
A core challenge across our linear combination layers is that correct terms may have relatively large coefficients, and the final mean absolute error for noiseless systems approaches 0. With this in mind, we note that we would ideally like the linear combination (and dense output) layers to function as feature selectors, which motivated us to implement $L_1$ regularization. We found, however, that the linear growth in error as coefficients increased was still somewhat unfavorable, so we instead implemented $L_{1/2}$ regularization~\cite{Xu2010}. $L_{1/2}$ regularization allows for greater sparsity than $L_1$ regularization while simultaneously providing diminishing penalties for increasing weight magnitudes. The regularization loss can be calculated as
\begin{equation}\label{l_half}
\alpha_1\sum_{i = 1}^{n+1}\left|w_i\right|^{\nicefrac{1}{2}} \equiv \alpha_1 L_{\nicefrac{1}{2}}.
\end{equation}
Here, $\alpha_1$ is a tuning parameter to modify the relative influence of this regularization. Linear combinations do not contribute significantly to the problem of exploding gradients in our application, and it is straightforward to initialize the system to small values; any standard initialization protocol with zero mean and non-constant initialization of weights is acceptable. We use hyperparameter value $\alpha_1 = 0.05$.

\subsubsection{Polynomial Combinations}
Unlike the linear combinations considered above, the polynomial combination layer does not correspond to sparsity in the output when all weights are 0. Rather, this corresponds to a constant value. Considering parsimony more broadly, we prefer simpler polynomial terms to more complicated ones: the lower the overall degree of the expression, the more we generally prefer it. For example, we view $x_1^2$ as preferable to both $x_1^5$ and $x_1^2 x_2$; we would also prefer $x_1 x_2$ to $x_1^2 x_2$ or $x_1^5$. Our primary mechanism of modulating the sparsity in the output expression is the constant value of 2 we previously appended to our state to assist with scaling. If the exponent on 2 is very negative (say, -10), then that term is virtually negligible. We therefore designed a novel regularizer to reward highly negative values of the exponent on 2 while penalizing large magnitudes (negative or positive) of the exponents for the non-constant input values. Formally, our regularizer loss is calculated as
\begin{equation}\label{l_poly}
\alpha_2 1.1^{\sum_{i = 1}^{n} \left|w_i\right| +  w_{n+1}} \equiv \alpha_2 L_{poly},
\end{equation}
which approaches 0 as $w_{n+1}\to-\infty$ and approaches infinity as the magnitudes of the non-constant values' exponents increase. We note that terms in the network model $\ftilde(\cdot)$ may be negligible in magnitude due to regularization but are not discarded unless deemed insignificant by the information criterion in the selected network model $\ftilde_A(\cdot)$, as described in Section~\ref{section_information}.

Unlike the case of the linear combination operational layers, here our initialization is critically important for ensuring stability of the system as the number of stacks increases. As noted previously, we generally aim for the initial value of each layer's output to be near-constant; in the case of the polynomial combination layer, this corresponds to a zero vector weight matrix, so we initialize the layer with a normal distribution about a mean of $0$. For deeply-stacked implementations, it is important that the standard deviation of this distribution is small in magnitude to avoid numerical issues. We use the hyperparameter value $\alpha_2 = 0.01$.

\subsubsection{Simple Products}

Regularization is less of a concern for the simple product operational layer because the output structure is already significantly constrained. All simple product layers output polynomials with degree no higher than the number of input terms, and the degree of a given variable within that expression is always either 1 or 0. As such, the enforcement of sparsity as it relates to these expressions can be safely applied at the dense output layer stage, rather than within the operational layer.

In contrast, initialization of the simple product operational layers is vital to the success of multi-stack implementations. Standard, mean-zero initialization (as was used elsewhere) will cause significant numerical issues in deeper networks because a weight of 0 does not correspond to a constant value, but rather to terms of the form $0.5+0.5x_1$. Instead, we center our initialization around a negative number, with the standard deviation and mean both decreasing as we wish to deepen our network. For systems with relatively few stacks, we found a mean of $-1.0$ to sufficiently protect against exploding errors and gradients.

\subsubsection{Common Operators}

The common operational layers are split into two sublayers of weights: first a sublayer performing a linear combination of the input expressions, then a sublayer performing a linear combination of the outputs of the operators acting on the first linear combination. Applying $L_{1/2}$ regularization to the first of these sublayers will help with sparsity and mitigate the risks of exploding gradients, but it is insufficient to promote parsimony since some operators possess nonzero values for zero-valued inputs. As such, we apply $L_{1/2}$ regularization to both the first and second sublayers of our system. However, we found in our hyperparameter studies that a reduced weight of $\alpha_3 = 0.0375$ keeps these functions from being penalized too much and thus we include it alongside Eg.~\ref{ops_eq} as
\begin{equation}\label{l_ops}
\alpha_3 \sum_{i = 1}^3\left|w_i\right|^{\nicefrac{1}{2}} \equiv \alpha_3 L_{ops}.
\end{equation}
Initialization is important here as well, and we again make sure to initialize about mean zero. However, we take an extra precaution here of also decreasing the standard deviation of our initialization distribution as we proceed through subsequent stacks of our network to further reduce the risk of initially high gradients and losses.

\subsection{Hyperparameters}

We use certain hyperparameter values (number of stacks and layers) to identify a variety of systems, which is analogous to using the same hyperparameters (width and depth) in a deep neural network to perform regression on entirely different datasets. As with traditional machine learning, setting the hyperparameters is not an exact science. Even when set after numerical experiments for one system such that the desired metrics are at their best, using data from a different system often requires re-tuning the hyperparameters. However, we found that a set of hyperparameters such as $K =1$ stack with $L = 10$ layers, or $K = 2$ stacks with $L = 2$ layers, works well across multiple systems, as demonstrated in subsequent sections. This is in conjunction with the same common operators of exponential, sinusoidal, and sign functions. Thus, we use the same hyperparameters to identify systems with fixed-points, limit cycles, and chaotic attractors.

When using $K = 1$ stack, we found $L < 10$ layers to be sufficient for some systems, but all systems were identifiable using $K=1$ stack and $L=10$ layers. This is surprising as many of the systems are, in principle, sufficiently expressed with just $L = 1$ or $L = 2$ layers. Such overparametrization is common in deep neural networks with hundreds of millions of parameters, but so is the issue of overfitting to training data. In contrast, our overparametrization seems to produce equations that are generalizable beyond the training data. This phenomenon of overparametrization beyond conventionally accepted number of parameters optimal to the trade-off between bias and variance~\cite{mohri2018foundations}, has recently been reported~\cite{double_descent} and is an active area of research.

The hyperparameters used in demonstrations of \mbox{SymANNTEx} presented in this work are summarized in Table \ref{hyperparameters}. These examples span various dynamic behaviors like hyperbolic (Takens-Bogdanov) and elliptic (simple pendulum) fixed points, limit cycles (FitzHugh-Nagumo, chemical kinetics) and a continuum of periodic orbits (simple pendulum), fast-slow dynamics (FitzHugh-Nagumo), and chaos (Lorenz, R\"{o}ssler, Chua), governed by equations with sinusoidal (simple pendulum), exponential (chemical kinetics), and polynomial terms. We see that such a wide array of systems are identified by SymANNTEx using nearly the same hyperparameters of training.

\begin{table*}[t]
\centering
\caption{Hyperparameters for demonstrated examples}\label{hyperparameters}
\begin{ruledtabular}
\begin{tabular}{cccccccc}
\multicolumn{1}{|c|}{Parameter} & Pendulum & Chua & Takens-Bogdanov & R\"{o}ssler & Lorenz & FitzHugh-Nagumo & \multicolumn{1}{c|}{Chemical Kinetics}\\
\hline
\multicolumn{1}{|c|}{Stacks} & \multicolumn{5}{c}{$1$} & \multicolumn{2}{|c|}{$2$} \\
\multicolumn{1}{|c|}{Layers} & \multicolumn{5}{c}{$10$} & \multicolumn{2}{|c|}{$1$} \\
\hline
\multicolumn{1}{|c|}{Learning rate constant} & \multicolumn{2}{c}{$0.032$} & \multicolumn{5}{|c|}{$0.01$}\\
\hline
\multicolumn{1}{|c|}{$L_{1/2}$ regularization weight} & \multicolumn{7}{c|}{$0.05$}\\
\multicolumn{1}{|c|}{$L_{poly}$ regularization weight} & \multicolumn{7}{c|}{$0.01$} \\
\multicolumn{1}{|c|}{$L_{ops}$ regularization weight} & \multicolumn{7}{c|}{$0.0375$} \\
\multicolumn{1}{|c|}{AIC tolerances} & \multicolumn{7}{c|}{Eq.~\ref{tolerance_intervals}} \\
\hline
\multicolumn{1}{|c|}{Number of data points} & \multicolumn{6}{c|}{$1000$} & \multicolumn{1}{c|}{$16000$} \\
\multicolumn{1}{|c|}{Training : test split} & \multicolumn{6}{c|}{$800:200$} & \multicolumn{1}{c|}{$12000:4000$} \\
\hline
\multicolumn{1}{|c|}{Training instances} & \multicolumn{7}{c|}{$5$ ($1$ per split)}\\
\hline
\multicolumn{1}{|c|}{Training Epochs} & \multicolumn{1}{c|}{100} & \multicolumn{1}{c}{6400} & \multicolumn{1}{|c|}{25} & \multicolumn{1}{c|}{400} & \multicolumn{1}{c}{6400} & \multicolumn{1}{|c|}{3200} & \multicolumn{1}{c|}{6400}\\ 
\hline
\multicolumn{1}{|c|}{Batch size} & \multicolumn{7}{c|}{$\min\{\sqrt{training \text{ } split}, \ 32\}$}\\
\hline
\multicolumn{1}{|c|}{Training Noise} & \multicolumn{6}{c|}{$\sigma_1 = 0.01$ (up-to 0.025 in some examples, $\sigma_2 = 0.01$ (up-to 0.075 in some examples)} & \multicolumn{1}{c|}{$\sigma_1 = \sigma_2 = 0.001$}\\
\end{tabular}
\end{ruledtabular}
\end{table*}

\subsection{Loss function and Training}

We use a combination of the Mean Absolute Error (MAE) and sublayer-dependent regularization as the loss which we seek to minimize using iterative optimizers, initialized with layer-dependent initial weights. While regression in Euclidean spaces is widely done using the Mean Squared Error (MSE) with $L_1$ regularization, we found (as shown in Sec.~\ref{comparison}) this conventional approach to be 1) lacking desired parsimony and 2) giving enormous errors aggravating the optimization. Instead, we use $L_{1/2}$ regularization~\cite{Xu2010} which is known to promote sparsity to a greater extent than $L_1$ does, but at the cost of losing convexity. 
The loss function is given by
\begin{equation}\label{loss}
        \text{Loss} = \text{MAE} + \sum_{k = 1}^K \sum_{l = 1}^L (\alpha_1 L_{1/2}^{(k,l)} + \alpha_2 L_{poly}^{(k,l)} + \alpha_3 L_{ops}^{(k,l)}),
\end{equation}
where
\[\text{MAE} = \frac{1}{m} \sum_{j = 1}^m \| \dot{\x}^{[j]}_{\sigma_2} - \ftilde(\x^{[j]}_{\sigma_1})\|_1,
\]
and $L_{1/2}^{(k,l)}, \ L_{poly}^{(k,l)}$ and $L_{ops}^{(k,l)}$ are custom regularizers at each sublayer (denoted by the superscript) in the stacks $k \in \{1,\cdots,K\}$ and operational layers $l \in \{ 1,\cdots,L\}$. $\alpha_1, \ \alpha_2$ and $\alpha_3$ are user-defined regularization weights. The \textit{possibly} longer computational time taken to minimize a non-convex loss is outweighed by the greater sparsity in fewer iterations of training. Iterations of training are commonly measured in ``epochs" where one epoch is when each of the data-points has been used once in the optimization. Greater sparsity in fewer epochs is desirable for increased interpretability.

We use the Adam optimizer~\cite{kingma2015adam} to minimize the loss function. We use a $4:1$ randomized split of the data for training and testing, respectively, in a \verb|Kfold| routine that uses different randomized initial weights for each of the 5 instances of the optimization. This 5-fold optimization increases the likelihood of an optimal set of network weights by concealing a different fifth of the data in each of the five training runs, and choosing that which has the smallest magnitude of the loss-function among all five of the instances. Although it is an algorithm with adaptive learning rate, it has a learning rate constant. We found the commonly used default learning rate constant of $0.001$ to be ``insufficient'' even with $10^5$ epochs yet a learning rate constant of $0.01$ seems to yield convergence well within $10^4$ epochs. This larger learning rate constant and the spurts of instantaneous increase in the loss function value at some iterations suggests that we could be observing what has recently been reported as the edge of stability~\cite{edge_of_stability}. Our investigations with learning rate constants of $\{0.001, \ 0.004, \ 0.016\}$ in our examples, each with datasets of $\{1000, \ 4000, \ 16000\}$ indicate the lack of a linear relationship between step-size and number of epochs taken for convergence to the expected equations. At slower learning rates, we found that the optimization seems to fall short of sparsity even when run for longer. In the interest of practical applications, we use higher learning rates which seem to remedy this over fewer training epochs. This is especially apparent with the Lorenz and FitzHugh-Nagumo models. While some models can be identified from as few as $25$ epochs, we train some examples for up to $6400$ epochs. Thus, we have used a learning rate constant of $0.01$ in almost all of our examples. It works for all the dataset sizes used in our investigation so we show results using the smallest datasets which have only $1000$ data-points. These hyperparameters are summarized in Table~\ref{hyperparameters}.

The learning rate for the simple pendulum has been increased by a factor of three. We found that the optimizer identifies a harmonic oscillator even when trained for $12800$ epochs with a learning rate of $0.01$, whereas increasing the learning rate to $0.032$ identifies Eq.~\ref{identified_pendulum} in as few as $25$ epochs. This indicates that the harmonic oscillator is a local minimum for the optimization in parameter space which is also the classical linear approximation to the simple pendulum. Reducing the regularization weight $\alpha_3$ (Eq.~\ref{l_ops}) of the common operators in Eq.~\ref{ops_eq} alleviates this but makes the regularization sub-optimal in identifying the other examples. We also found that the higher learning rate worked better for the Chua double scroll system.

\section{Information criterion for parsimonious models}
\label{section_information}
We seek parsimonious expressions for the network model for given (noisy) data, because this increases interpretability. 
In practice, we find that the network model has coefficients that 1) do not exactly match those of the ground truth model, although they are often approximately equal, and 2) there are a large number of small but non-zero coefficients despite sparsity-promoting regularization of the empirical loss. Thus, we need a ``post-processing'' step that tunes the coefficients of the network model, including the possibility of setting coefficients of terms to zero. 
This can be posed as a problem in model selection.

We use the Akaike Information Criterion (AIC)~\cite{akaike1974new} to fine-tune the network model's parameters. 
If an approximate model can be expressed using fewer terms in the equation, it could be viewed as carrying similar informativity but with fewer parameters. Information theoretic criteria have been used in statistics for model selection given a possible set of models, including with SINDy~\cite{mangan2017model}; in particular, AIC can be viewed as quantifying the complexity of a model via the number of terms in its equation alongside the residual relative to the actual data, found from maximum likelihood estimates. As we work in a Euclidean space, we use the MSE for $m$ datapoints. AIC uses the log-likelihood estimate alongside a correction term that accounts for finite data size~\cite{model_selection}
\begin{equation}\label{AIC}
        \text{AIC} =  \big(2 P + m\log_e(\text{MSE})\big) + 2\frac{(P+1)(P+2)}{m - P - 2},
\end{equation}
where
\[\text{MSE} = \frac{1}{m} \sum_{j = 1}^m \| \dot{\x}^{[j]}_{\sigma_2} - \ftilde_S(\x^{[j]}_{\sigma_1})\|^2_2,
\]
and $P$ is the number of estimated parameters.
Here, $\ftilde_S(\cdot)$ is a ``simpler'' model obtained by using the Python library SymPy's \verb|nsimplify| that rounds parameters of the network-model $\ftilde(\cdot)$ to within a given tolerance. This command also searches for commonly occurring constants like $\pi, e,$ etc. to within the tolerance. We compute the AIC scores for the network-model over a range of user-specified tolerances, and select the model $\ftilde_A(\cdot)$ with the lowest AIC score. The tolerances considered are $11$ logarithmically spaced numbers between $0.01$ and $1$ as
\begin{multline}\label{tolerance_intervals}
    \text{Tolerances used} \equiv [0.01, \ 0.0158, \ 0.0251, \ 0.0398, \\
    \ 0.0631, \ 0.1, \ 0.1585, \ 0.2512, \ 0.3981, \ 0.6310, \ 1].
\end{multline}
Rounding parameters of the network-model $\ftilde(\cdot)$ to each of these 11 tolerances could give up to 11 models differing in parameters. The AIC score is computed on each of these models and the one with the lowest score is the selected network model $\ftilde_A(\cdot)$. The effect of the above rounding-off can sometimes be drastic when used in conjunction with sparsity inducing regularization because many of the network model parameters could be close to 0. By following the aforementioned procedure on a network model $\ftilde(\cdot)$ with $322$ parameters, the Lorenz system is identified exactly as $\ftilde_A(\cdot)$ with 5 parameters, or the simple pendulum is identified closely with 2 parameters.

\section{Examples}

We have successfully identified a number of models whose characteristics span a variety of dynamic behaviors of scientific and engineering importance, using noisy data. Here we provide some of those examples to demonstrate identification of systems with stable, unstable, and saddle fixed points, a limit cycle, and chaotic behavior. Some of these systems also show time-scale separation via fast initial transient dynamics onto a slow manifold. We also compare trajectories generated from a random initial condition by the selected network models to that of the ground truth model.

In this paper we present examples trained on data with $1\%$ noise (as defined in Eq.~\ref{noisy_data}) but we first show how that added noise translates to the metrics of performance used. Although we use the MAE in the loss function (Eq.~\ref{loss}) that we optimize for, our data lies in a Euclidean space so we look at the Root Mean Squared Error (RMSE) which is the $l^2$ error averaged over our dataset. We already defined MSE in Eq.~\ref{AIC} and the RMSE is the square-root of the same but it can vary across different dynamical systems, so we consider the relative-RMSE which is defined as
\begin{equation}\label{rmse}
    \text{Relative-RMSE} = \sqrt{\frac{\sum_{j = 1}^m \| \ftilde_A(\x^{[j]}_{\sigma_1}) - \dot{\x}^{[j]}_{\sigma_2} \|_2^2}{\sum_{j = 1}^m  \| \dot{\x}^{[j]}_{\sigma_2} \|_2^2}}.
\end{equation}
The above would be non-zero even if the ground truth model is identified exactly, since $E(\dot{\x}_{\sigma_2}) \neq E(\f(\x_{\sigma_1}))$ from Sec.~\ref{prelims}. While in principle it is possible that some $\ftilde_A(\cdot)$ could over-fit any finite dataset, we know that the least relative-RMSE of interest is obtained if $\ftilde_A \equiv \f$. This can be used to lower-bound the relative-RMSE. Computing it theoretically becomes cumbersome, so we instead provide numerical estimates in Table~\ref{rmse_table} using Eq.~\ref{rmse} for each example, in the limit of large data, to better approximate the errors.
\begin{table}[t]
\centering
\caption{Approximate relative-RMSE for each of the examples computed using ground truth models over $100,000$ data-points spread across $1000$ trajectories for the non-chaotic examples and over single trajectories computed for $T = 1000$ time-units for the chaotic examples}\label{rmse_table}
\begin{ruledtabular}
\begin{tabular}{cc}
Model & $\approx$ Relative-RMSE\\
\hline
Takens-Bogdanov &  $2.04\%$\\
Simple Pendulum &  $4.37\%$\\
R\"{o}ssler &  $2.69\%$\\
Lorenz &  $2.67\%$\\
FitzHugh-Nagumo &  $2.64\%$\\
Chemical kinetics & $0.44\%$\\
Chua double-scroll & $1.83\%$\\
\end{tabular}
\end{ruledtabular}
\end{table}

\subsection{Takens-Bogdanov normal form}

We consider the Takens-Bogdanov normal form~\cite{guckenheimer2013nonlinear} to demonstrate identification of a model using data for a system that has fixed points with different stability characteristics
\begin{subequations}
\begin{align}
    \dot{x} &= y, \label{tb_x}\\
    \dot{y} &= \mu_1 + \mu_2 \; y + x^2 + x \; y. \label{tb_y}
\end{align}
\end{subequations}
This is the prototypical example of a dynamical system which is close to parameter values for which a fixed point has a double zero-eigenvalue, and possible special solutions can include fixed points $P_\pm \equiv (\pm \sqrt{- \mu_1}, 0)$ and limit cycles. 
We consider the parameter values $\mu_1 = -4.41, \ \mu_2 = 1.5,$ which is a first-order approximation to values for which there exists a homoclinic orbit~\cite{guckenheimer2013nonlinear}. For the considered parameter values there is a stable fixed point $P_-$, an unstable (saddle) fixed point, no limit cycles, and trajectories with initial conditions outside of the basin of attraction of $P_-$ grow without bound.

\begin{figure}[htbp]
\centering
\hspace{-14.5pt}\includegraphics[width=0.855\linewidth]{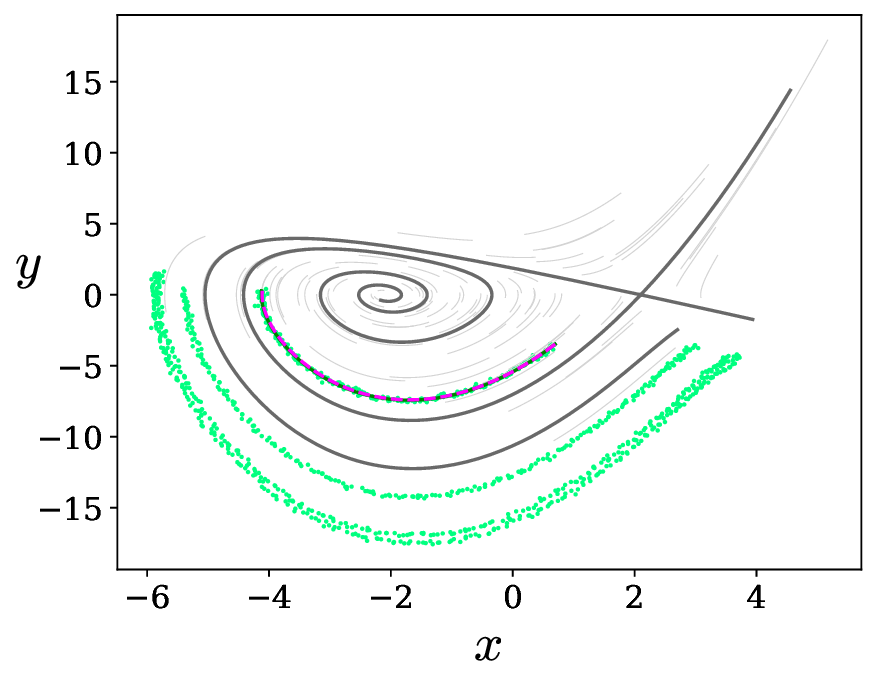}
\newline
(a) Phase portrait
\includegraphics[width=0.905\linewidth]{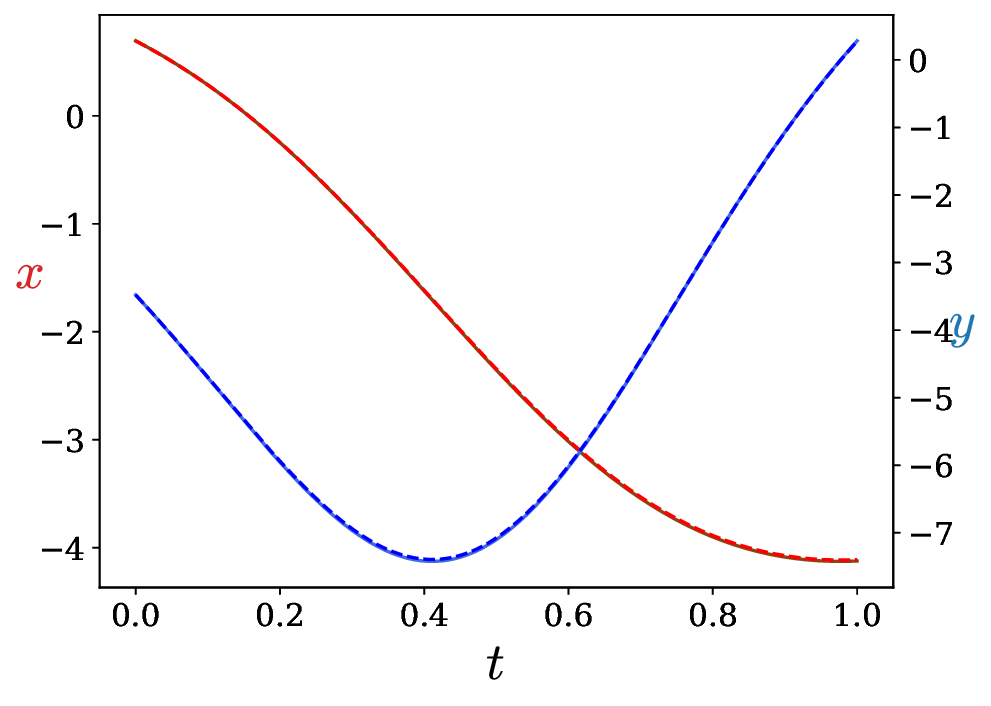}
\newline
(b) Time series
\caption{(a) State-space comparison of the trajectory generated (dashed magenta) by Eqs.~(\ref{tb_main_x}-\ref{tb_main_y}) estimated at $\X$ (light green) with that generated by the ground truth model (Eqs.~\ref{tb_x}-\ref{tb_y}) (solid green) which here is visually indistinguishable. The stable and unstable manifolds of the saddle fixed point are shown as dark gray lines, and other trajectories are shown as light gray lines for better depiction of the dynamics. (b) The same trajectory (dashed) plotted as a time series compared with that generated by the ground truth model (solid), which are nearly identical.}
\label{tb_fig}
\end{figure}

To test our system identification algorithm, we consider as input the $4$ randomly initialized trajectories shown in Fig.~\ref{tb_fig}(a). There are a total of $1000$ data points ($250$ for each trajectory of length $T = 1$). We find that SymANNTEx identifies the model
\begin{subequations}
\begin{align}
    \dot{x} &= y,\label{tb_main_x}\\
    \dot{y} &= -4.333 + 1.5 \; y + x^2 + x \; y, \label{tb_main_y}
\end{align}
\end{subequations}
using the following network structures: 1) $L = 10$ layers with $K = 1$ stack with $236$ trainable parameters, and 2) $L = 1$ layer and $K = 2$ stacks with $68$ trainable parameters. This model is estimated after $25$ epochs on $\Y$ at $\X$ shown in Fig.~\ref{tb_fig} with a relative-RMSE of $\approx~2.10\%$. For reference, the ground truth model gives an error of $\approx~2.04\%$ on the dataset referred to in Table~\ref{rmse_table}. SymANNTEx estimated similar models up to noise levels of $\sigma_1 = 0.02, \sigma_2 = 0.075$.

\subsection{Simple pendulum}

To demonstrate the ability to identify models with sinusoidal terms, we consider the simple pendulum
\begin{equation}\label{pendulum_eqn}
    \ddot{\theta} = -\frac{g}{l} \sin \theta,
\end{equation}
with $g = 9.81, \ l = 2$. As our algorithm is designed to approximate first-order ordinary differential equations (ODEs) as functions of state, we reformulate the above second-order ODE equation as two first-order ODEs, which is a classical technique in analyzing dynamical systems~\cite{guckenheimer2013nonlinear}. The simple pendulum has a continuum of time-periods ranging from $2\pi\sqrt{\nicefrac{l}{g}}$ near the fixed-point at $(0,0)$ to infinity near the unstable fixed-points $(\pm k\pi, 0)$. We remark that this continuum of time periods is reflected in the Koopman operator perspective as continuous spectrum~\cite{mezic2020spectrum} which makes analysis using DMD challenging~\cite{lusch2018deep}.

\begin{figure}[htbp]
\centering
\hspace{-11pt}\includegraphics[width=0.835\linewidth]{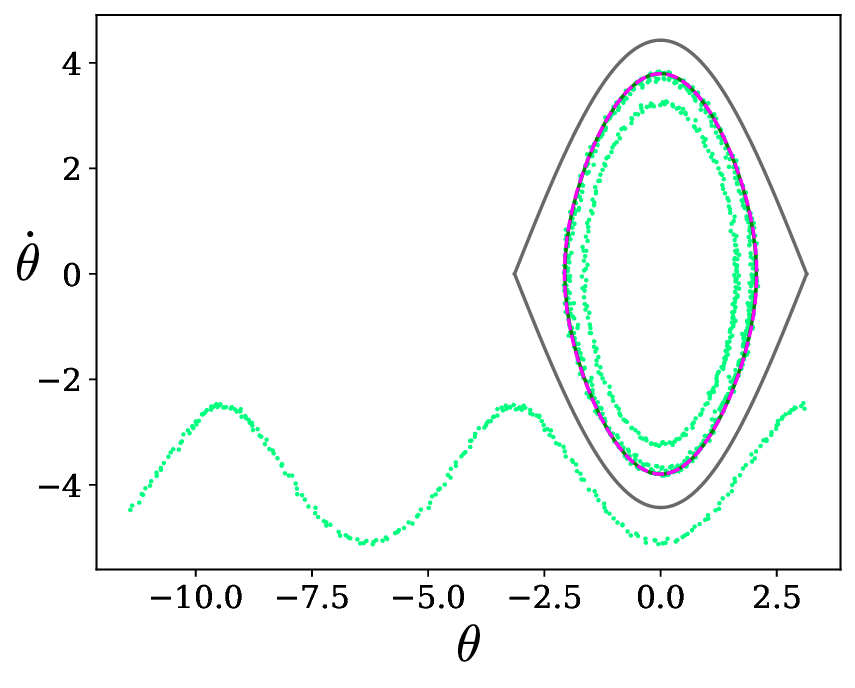}
\newline
(a) Phase portrait
\includegraphics[width=0.9\linewidth]{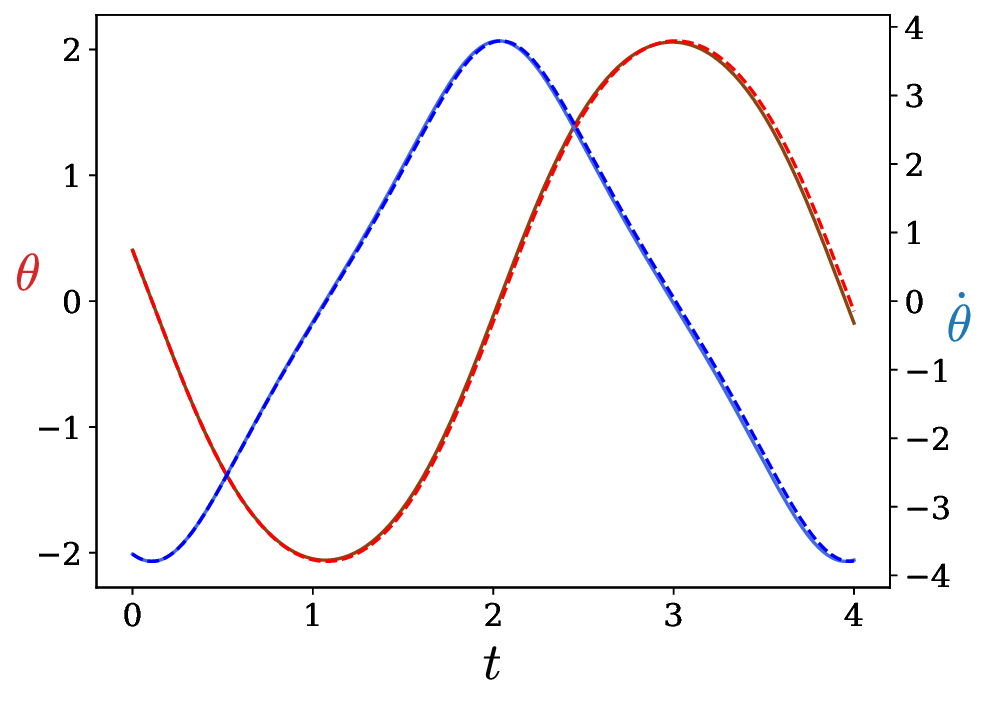}
\newline
(b) Time series
\caption{(a) State-space comparison of the trajectory generated (dashed magenta) by Eqs.~\ref{identified_pendulum} estimated at $\X$ (light green) with that generated by the ground truth model (Eq.~\ref{pendulum_eqn}) (solid green), which here is visually indistinguishable. The separatrices between angular displacement $-\pi$ and $\pi$ are shown as dark grey lines. (b) The same trajectory (dashed) plotted as a time series compared with that generated by the ground truth model (solid), which are nearly identical.}
\label{pendulum_fig}
\end{figure}

To test our algorithm, we consider as input the $4$ randomly initialized trajectories shown in Fig~\ref{pendulum_fig}. There are a total of $1000$ data points ($250$ for each trajectory of length $T = 4$). We find that our algorithm identifies the model in the state-space form, denoting $x \equiv \theta$ and $y \equiv \dot{\theta}$, as
\begin{subequations}\label{identified_pendulum}
\begin{align}
    \dot{x} &= y\\
    \dot{y} &= -4.857 \sin x,
\end{align}
\end{subequations}
using $L = 10$ layers and $K = 1$ stack with $236$ trainable parameters after $100$ epochs on training data $\Y$ at the input data $\X$ with a relative-RMSE of $\approx~4.45\%$. For reference, the ground truth model gives an error of $\approx~4.37\%$ on the dataset referred to in Table~\ref{rmse_table}. Here, we remark that a higher learning rate of $0.032$ is used to identify the above equations in as few as $25$ epochs as compared to a a learning rate of $0.01$ only giving a linear approximation even after $12800$ epochs indicating a large region of local minima. SymANNTEx estimated similar models up to noise levels of $\sigma_1 = 0.02, \sigma_2 = 0.075$.

\subsection{R\"{o}ssler equations}

We now demonstrate the identification of chaotic dynamical systems from noisy data. Identifying the flow-map or solution of state of chaotic systems is challenging due to their sensitive dependence on initial conditions and usual lack of closed form solutions. Their analysis using DMD methods is also challenging due to the lack of discrete modes\cite{budi12}. Here we consider the R\"{o}ssler equations which find relevance in chaotic behavior of chemical reactions~\cite{rossler1976chaotic}
\begin{subequations}
\begin{align}
  \dot{x} &= -y - z, \label{rossler_x}\\
    \dot{y} &= x + 0.5 \; y,\\
    \dot{z} &= 2 + z \; (x - 4). \label{rossler_z}
\end{align}
\end{subequations}
This model has a chaotic attractor with sensitive dependence on initial conditions. It also has a region of state-space beyond the attractor's region of attraction where trajectories grow unboundedly.

\begin{figure}[htbp]
\centering
\hspace{-9.5pt}\includegraphics[width=0.84\linewidth]{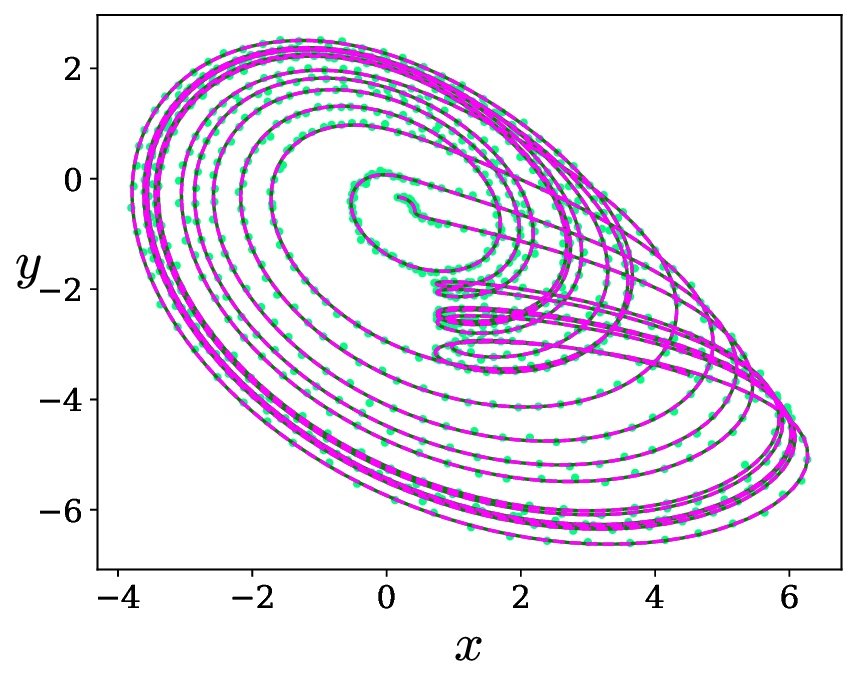}
\newline
(a) Phase portrait
\includegraphics[width=0.9\linewidth]{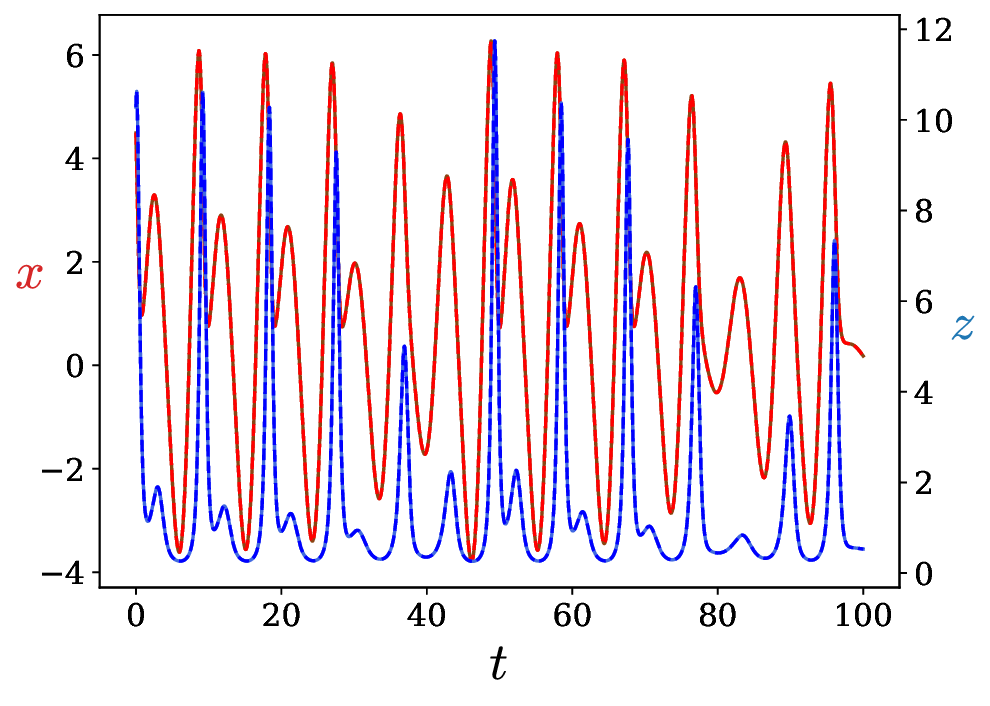}
\newline
(b) Time series
\caption{(a) State-space visualization of the trajectory generated by Eqs.~(\ref{rossler_x}-\ref{rossler_z}) (solid green). It coincides with that generated (dashed magenta) by the model identified at $\X$ (light green) as it is the same as the ground truth model. (b) The same trajectory (dashed) plotted as a time series coinciding with the ground truth. The time-evolution of the states is without any coherence due to the chaotic nature of the system. ($z$ is not shown in (a) and $y$ is not shown in (b) to avoid cluttering the figures).}
\label{rossler_fig}
\end{figure}

To test our algorithm, we consider $1000$ data points from a single randomly initialized trajectory (equally spaced in time on the time interval $[0, 100]$), as shown in Fig.~\ref{rossler_fig}. Since the terms in (\ref{rossler_x}-\ref{rossler_z}) can be represented without stacking, we use $L = 10$ layers and $K = 1$ stack with $322$ trainable parameters. This network was found to successfully identify the equations exactly, including the exact parameters, in $800$ epochs on training data $\Y$ at the input data $\X$ which has a relative-RMSE of $\approx~2.60\%$. We remark that the noise levels in Table~\ref{rmse_table} are estimates that vary slightly with the number of data points and are therefore not exact. SymANNTEx estimated exact model up to noise levels of $\sigma_1 = 0.01, \sigma_2 = 0.05$.

\subsection{Lorenz equations}

We further demonstrate the ability of SymANNTEx to identify chaotic dynamical systems from noisy data. We consider the classical Lorenz equations with standard parameters, originally used to model two-dimensional atmospheric convection~\cite{lorenz1963deterministic}
\begin{subequations}
\begin{align}
  \dot{x} &= 10 \; (y - x), \label{lorenz_x}\\
    \dot{y} &= x \; (28 - z) - y, \label{lorenz_y}\\
    \dot{z} &= x \; y - 8 z / 3. \label{lorenz_z}
\end{align}
\end{subequations}
This model has a chaotic attractor with sensitive dependence on initial conditions and a lack of closed form expression for its solution of state which manifests itself in forecasting time series. Reservoir computing\cite{bollt2021explaining} does an excellent job in predicting the state for many Lyapunov time-units but eventually diverges. Its associated Koopman operator lacks discrete spectrum\cite{arba17} and thus analysis using DMD is challenging with results continuously varying with the order of approximation.

\begin{figure}[htbp]
\centering
\hspace{-7.5pt}\includegraphics[width=0.855\linewidth]{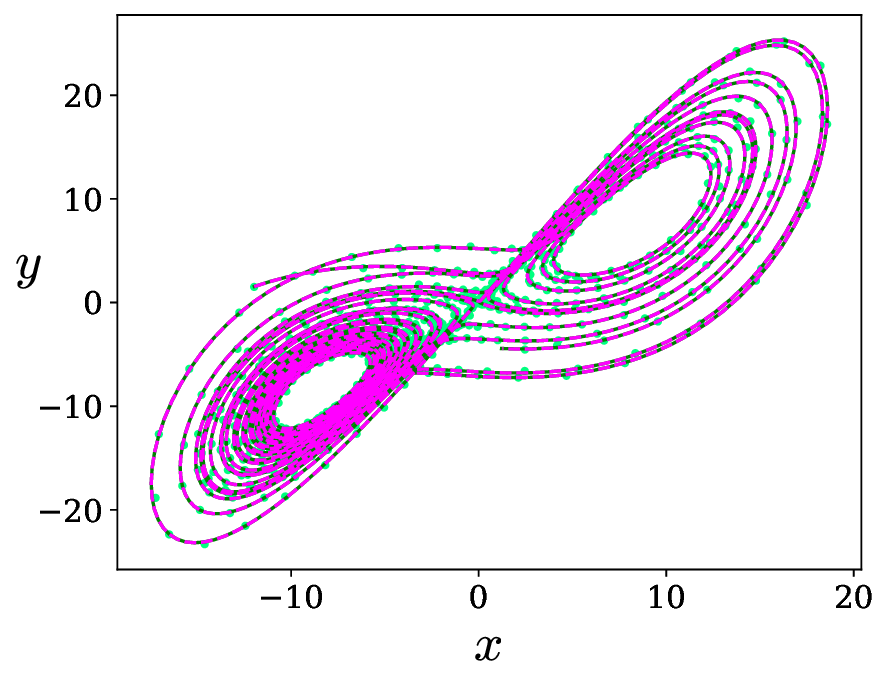}
\newline
(a) Phase portrait
\includegraphics[width=0.9\linewidth]{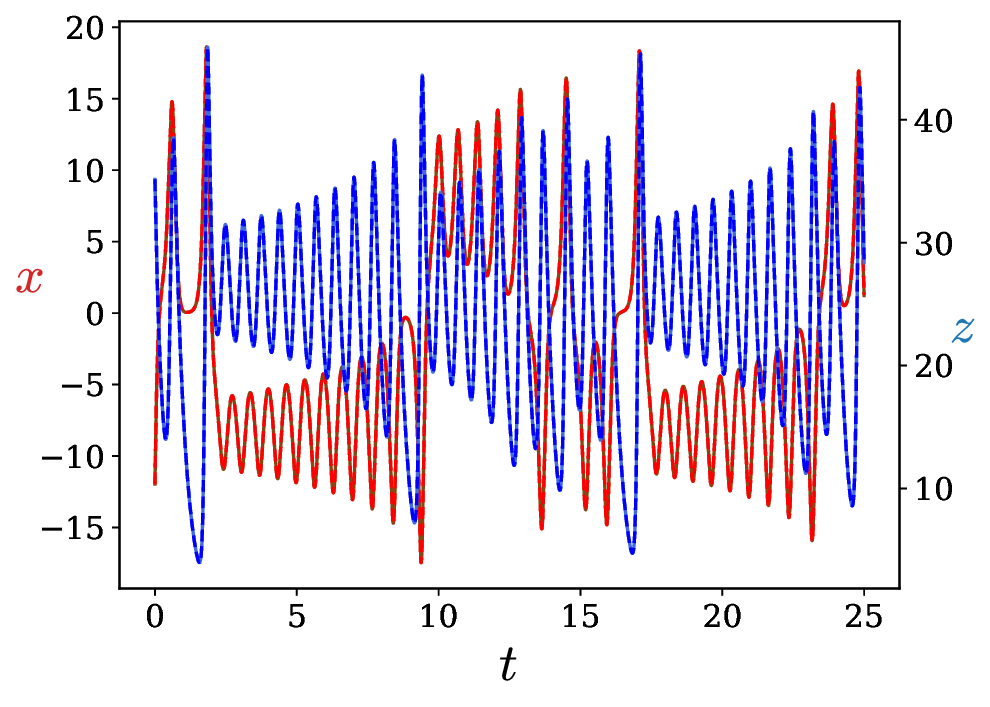}
\newline
(b) Time series
\caption{(a) State-space visualization of the trajectory generated by Eqs.~(\ref{lorenz_x}-\ref{lorenz_z}) (solid green). It coincides with that generated (dashed magenta) by the model identified at $\X$ (light green) as it is the same as the ground truth model. (b) The same trajectory (dashed) plotted as a time series coinciding with the ground truth. The time-evolution of the states is without any coherence due to the chaotic nature of the system. ($z$ is not shown in (a) and $y$ is not shown in (b) to avoid further cluttering the figures).}
\label{lorenz_fig}
\end{figure}

To test our algorithm, we consider $1000$ data points from a single randomly initialized trajectory (equally spaced in time on the time interval $[0, 25]$), as shown in Fig.~\ref{lorenz_fig}. Since the terms in (\ref{lorenz_x}-\ref{lorenz_z}) can be represented without stacking, we use $L = 10$ layers and $K = 1$ stack with $322$ trainable parameters. This network was found to successfully identify the equations exactly, including the exact parameters, in $6400$ epochs on training data $\Y$ at the input data $\X$ which has a relative-RMSE of $\approx~2.59\%$. We remark that the noise levels in Table~\ref{rmse_table} are estimates that vary slightly with the number of data points and are therefore not exact. This example with coefficients as large as $28$ alongside most of the other coefficients from the trainable parameters that get approximated to $0$ shows the wide range of parameters that can be identified across orders of magnitude. SymANNTEx estimated the exact model up to noise levels of $\sigma_1 = 0.01, \sigma_2 = 0.05$.

\subsection{FitzHugh-Nagumo equations}

We consider the FitzHugh-Nagumo model~\cite{fitzhugh1961impulses}
\begin{subequations}
\begin{align}
    \dot{v} &= v - v^3/3  - w + I, \label{FHN_v}\\
    \dot{w} &= \epsilon ( v + a w + b), \label{FHN_w}
\end{align}
\end{subequations}
to demonstrate identification of a system with oscillatory behavior and a separation of time-scales. This is a model for neural dynamics which shares key features with conductance-based models such as the Hodgkin-Huxley equations. Here we use the parameter values $I = 0.328, \ \epsilon = 0.08, \ a = -0.8, \ b = 0.7,$ for which there is a stable limit cycle with spiking behavior.

\begin{figure}[htbp]
\centering
\hspace{-23pt}\includegraphics[width=0.86\linewidth]{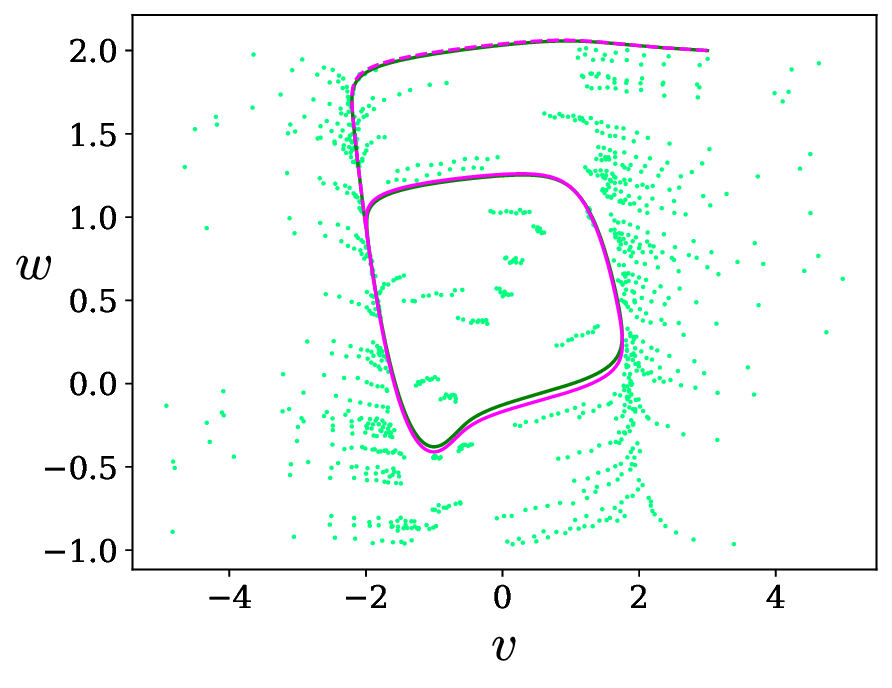}
\newline
(a) Phase portrait
\includegraphics[width=0.93\linewidth]{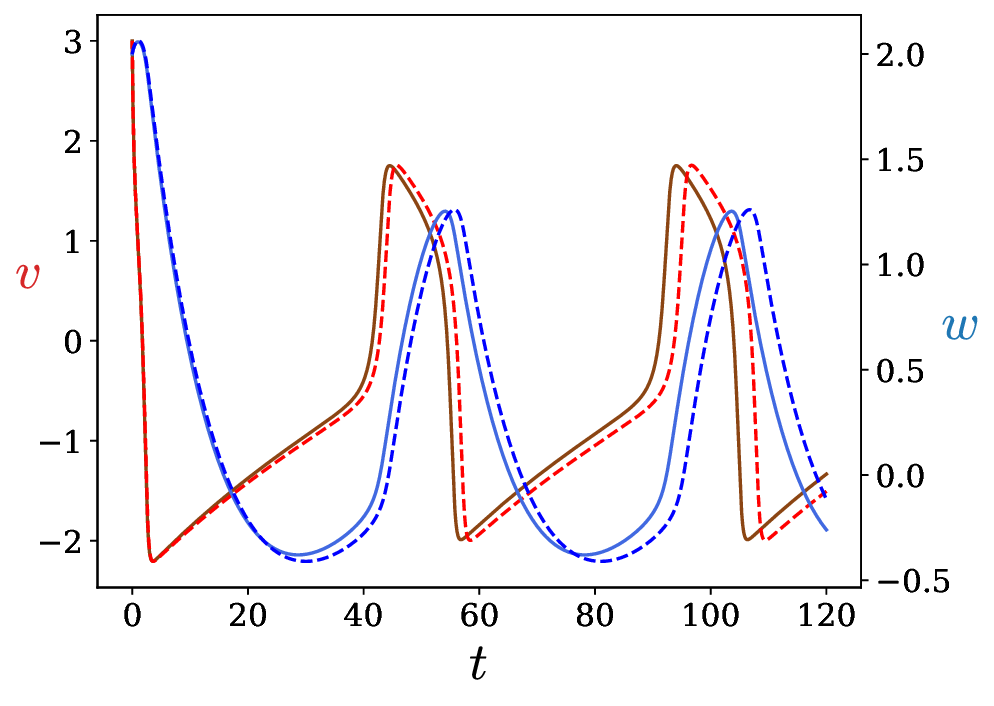}
\newline
(b) Time series
\caption{(a) State-space comparison of trajectory generated (dashed magenta) by Eqs.~(\ref{v_main}-\ref{w_main}) estimated at $\X$ (light green) with that generated by the ground truth model (Eqs.~\ref{FHN_v}-\ref{FHN_w}) (solid green). The fast initial transient of the trajectory and subsequently closely following the slow manifold onto the periodic orbit are almost exact. The periodic orbit itself is also reproduced closely. (b) Evolution of the same trajectories as a time series shows how the transients are almost exact and the time-period of the identified model (dashed) is slightly different from the ground truth model (solid).}
\label{fhn_fig}
\end{figure}

To test our algorithm, we consider as input $100$ randomly initialized trajectories shown in Fig.~\ref{fhn_fig}(a).  There are a total of $1000$ data points ($10$ for each trajectory of length $T = 1$). We find that SymANNTEx identifies the model
\begin{subequations}
\begin{align}
    \dot{v} &= 0.998 \ v - v^3/3  - 0.971 \ w + 0.311 - 0.006 \ v^2,\label{v_main}\\
    \dot{w} &= 0.078 \ v - 0.06 \ w + 0.054\label{w_main}
\end{align}
\end{subequations}
using $L = 1$ layer and $K = 2$ stacks with $68$ trainable parameters. We note that stacking is necessary because powers of state are taken on the absolute value, so the cubic term in (Eq.~\ref{FHN_v}) cannot be exactly formed using a single stack: $v^3 = |v|^2 \times v$. This model is identified in $3200$ epochs on training data $\Y$ at the input data points $\X$ with a relative-RMSE of $\approx~2.40\%$. For reference, the ground truth model gives an error of $\approx~2.64\%$ on the dataset referred to in Table~\ref{rmse_table} indicating a \textit{slight} overfit. Yet, from Fig.~\ref{fhn_fig}, we remark that the phase portraits nearly coincides with the ground truth model including the dynamic range of oscillations and the fast and slow transients. In particular, the time-period of the limit cycle behavior is also nearly the same, as seen from Fig.~\ref{fhn_fig}(b). SymANNTEx estimated similar models up to noise levels of $\sigma_1 = 0.025, \sigma_2 = 0.025$.

\subsection{Chemical kinetics with Arrhenius rate dependence}

We consider the exponential approximation to Arrhenius rate law in chemical kinetics in the context of two first-order reactions~\cite{gray1990chemical}, which gives the nonlinear dependence of reaction rate on the temperature. When the rise in temperature in an exothermic reaction is small in comparison to the ambient temperature, the rate can be approximated as exponentially dependent on the rise in temperature. After scaling the governing equations of mass and energy by relevant quantities,
\begin{subequations}
\begin{align}
    \dot{\alpha} &= - \kappa \; \alpha \; \text{e}^{\theta} + \mu, \label{arr_alpha}\\
    \dot{\theta} &= \alpha \; \text{e}^{\theta} - \theta \label{arr_theta}
\end{align}
\end{subequations}
is the set of dimensionless equations, where $\alpha$ is the intermediate chemical concentration, and $\theta$ represents temperature rise. We consider reactant concentration measure $\mu = 0.1$ and reaction rate constant $\kappa = 0.07$ for which this system of equations has a stable limit cycle.

\begin{figure}[htbp]
\centering
\hspace{-6.5pt}\includegraphics[width=0.825\linewidth]{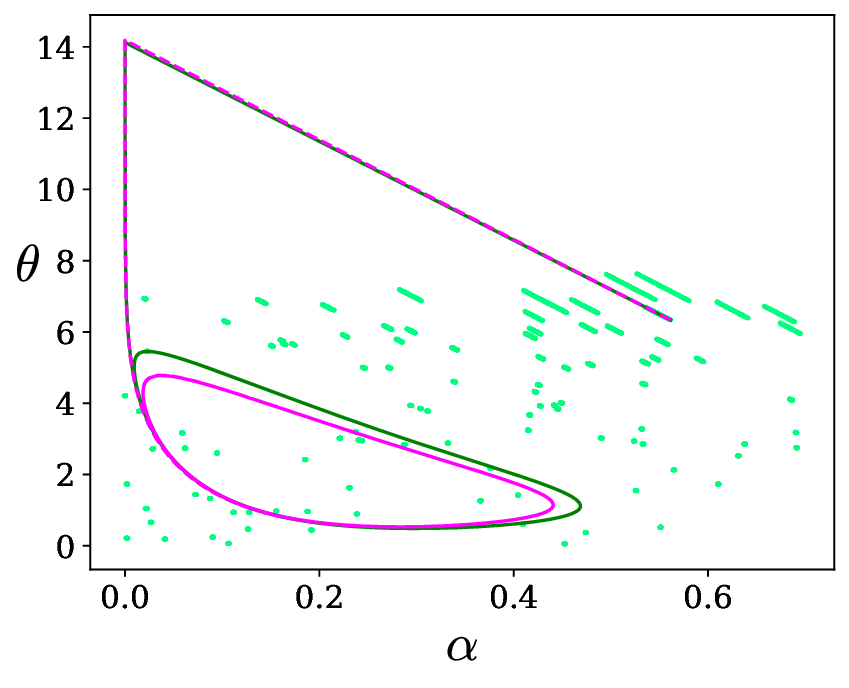}
\newline
(a) Phase portrait
\includegraphics[width=0.9\linewidth]{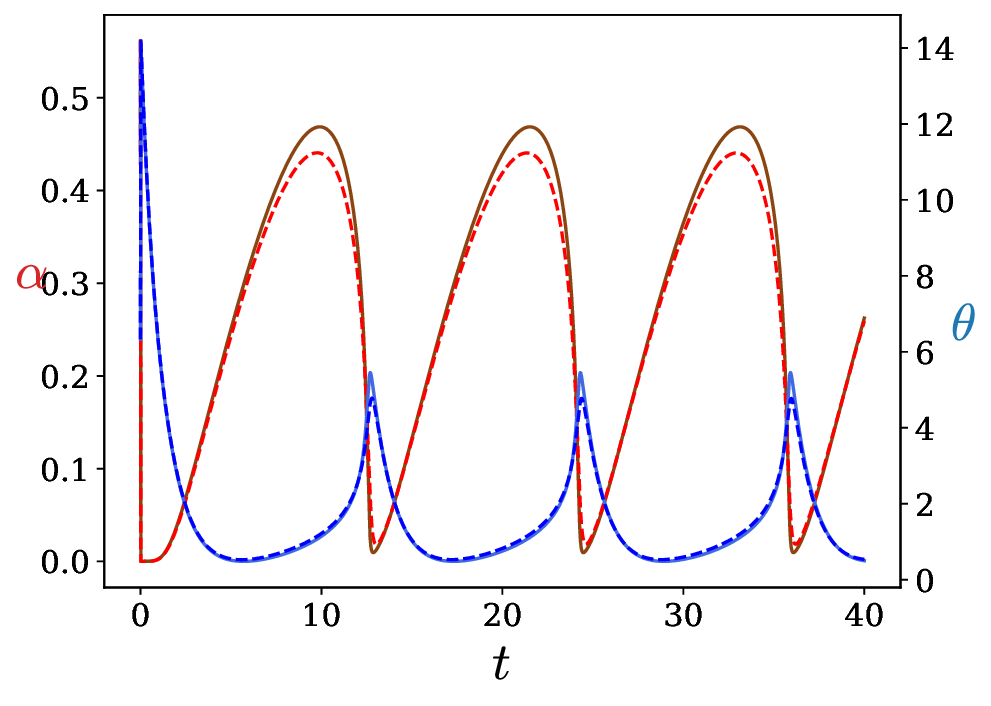}
\newline
(b) Time series
\caption{(a) State-space comparison of trajectory generated (dashed magenta) by Eqs.~(\ref{alpha_main}-\ref{theta_main}) estimated at $\X$ (light green) with that generated by the ground truth model (Eqs.~\ref{arr_alpha}-\ref{arr_theta}) (solid green). The initial transient of the trajectory onto the periodic orbit happens very quickly and is reproduced almost exactly. The periodic orbit itself is similar, capturing that the dynamic range in the two states differs nearly by an order of magnitude. (b) Evolution of the same trajectories in time shows how the transient and the time-period of the identified model (dashed) is almost exactly the same as that of the ground truth model (solid).}
\label{arr_fig}
\end{figure}

To test our algorithm, we consider as input $100$ trajectories shown in Fig.~\ref{arr_fig}(a). There are a total of $16000$ data points ($160$ for each trajectory of length $T = 0.001$). We find that SymANNTEx identifies

\begin{subequations}
\begin{align}
    \dot{\alpha} &= -0.07 \ \alpha \ \text{e}^{0.9 \ \theta} \; \frac{1.014 \ + \ 0.171 \ \theta}{\theta^{0.1} \ \text{e}^{0.1 \ \alpha}} \ + \ 0.1,\label{alpha_main}\\
    \dot{\theta} &= \alpha \ \text{e}^{0.9 \ \theta} \ \frac{1.028 \ + \ 0.175 \ \theta \ + \ 0.007 \ \theta^2}{\theta^{0.1} \ \text{e}^{0.1 \ \alpha}} \ - \ \theta,\label{theta_main}
\end{align}
\end{subequations}
using $L = 1$ layer and $K = 2$ stacks with $68$ trainable parameters. We note from Eqs.~(\ref{arr_alpha}-\ref{arr_theta}) that stacking is necessary because the exponential term $\text{e}^\theta$ is multiplied with a state $\alpha$ which cannot be exactly formed using a single stack. This model is identified in $6400$ epochs on training data $\Y$ at the input data points $\X$ with a relative-RMSE of $\approx~1.56\%$. For reference, the ground truth model gives an error of $\approx~0.44\%$ in Table~\ref{rmse_table}. The differences in Eqs.~(\ref{alpha_main}-\ref{theta_main}) compared to Eqs.~(\ref{arr_alpha}-\ref{arr_theta}) are terms with coefficients that are an order of magnitude smaller, which are even smaller when the range of $\X$ ($\alpha \in (0,0.7], \ \theta \in (0, 7]$) is considered.
We see from Fig.~\ref{arr_fig}(a) that the phase portraits are very similar with the ground truth model, including the dynamic range of oscillations. Although there is a slight mismatch in the fast transients, it is a range of state-space outside of $\X$ and on a faster time-scale than $\Y$. From Fig.~\ref{arr_fig}(b) we can also see that the time-period and dynamic range of oscillations is also nearly reproduced. We remark that the fast and slow transients are also reproduced closely. This also shows our algorithm's ability to identify datasets with large separation of time-scales. SymANNTEx estimated similar models up to noise levels of $\sigma_1 = 0.001, \sigma_2 = 0.01$.

\begin{figure}[htbp]
\centering
\hspace{-6.5pt}\includegraphics[width=0.825\linewidth]{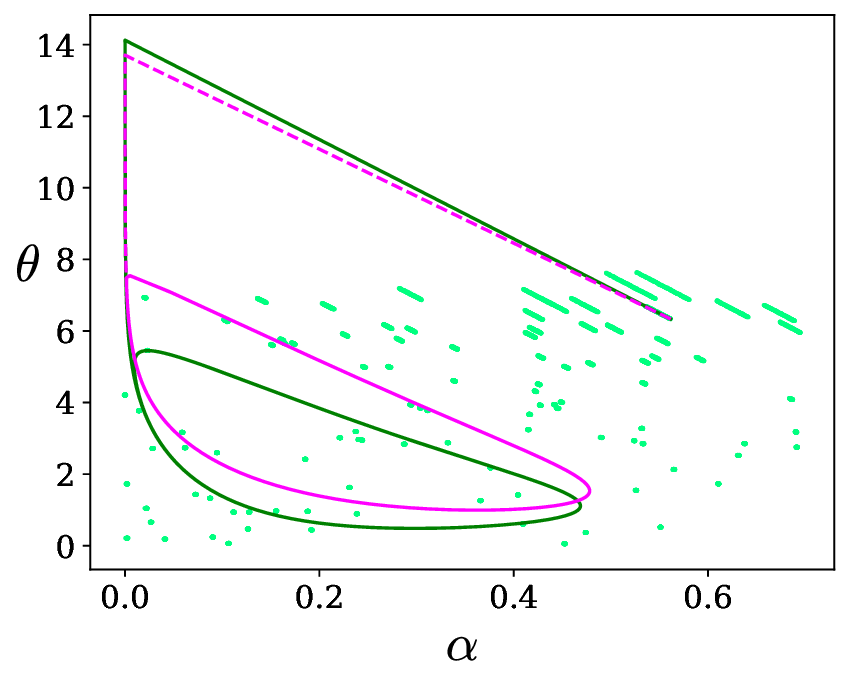}
\caption{State-space trajectory generated by Eqs.~(\ref{alpha_main_2}-\ref{theta_main_2}) (dashed magenta) compared to that generated by Eqs.~(\ref{arr_alpha}-\ref{arr_theta}) (solid green). Although the equations are visually similar, Eqs.~(\ref{alpha_main}-\ref{theta_main}) are a better approximation to the ground truth model in terms of the RMSE over $\X$ and the trajectories generated.}
\label{arr_fig_2}
\end{figure}

Below we also show another model identified for the same dataset but with a worse AIC score
\begin{subequations}
\begin{align}
    \dot{\alpha} &= -0.075 \ \alpha \ \text{e}^{\theta} + \ 0.167,\label{alpha_main_2}\\
    \dot{\theta} &= 1.016 \ \alpha \ \text{e}^{\theta} - \ \theta.\label{theta_main_2}
\end{align}
\end{subequations}
At first glance, these equations look more similar to the ground truth model but in comparison to Eqs.~\ref{alpha_main}-\ref{theta_main}, they have larger Relative-RMSE and, as seen from Fig.~\ref{arr_fig_2}, the trajectory they generate is not as close to that of the ground truth model. 

\subsection{Chua double-scroll attractor}
\begin{figure}
\centering
\hspace{-17pt}\includegraphics[width=0.86\linewidth]{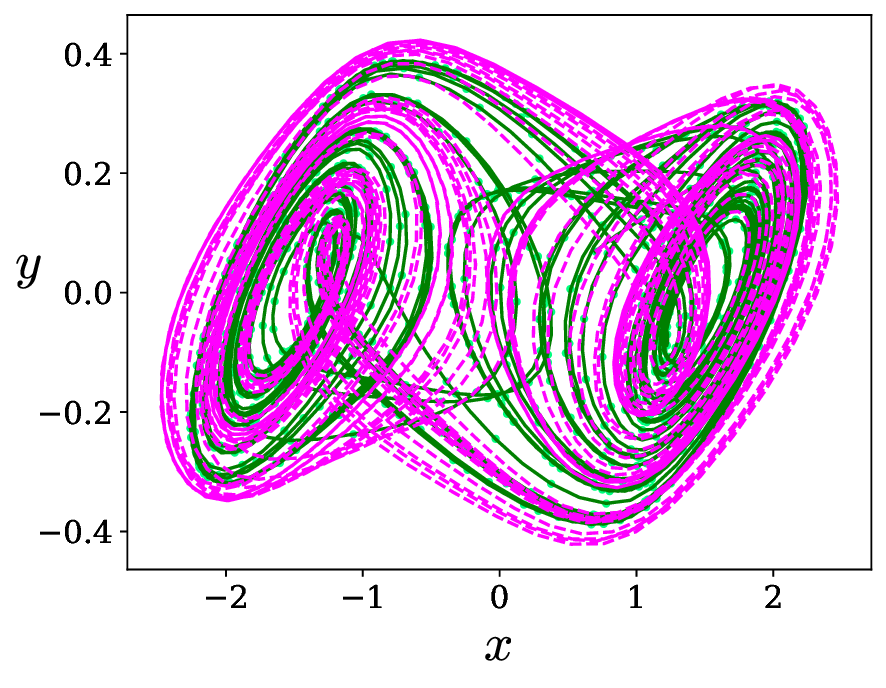}
\newline
(a) Phase portrait
\hspace{7.5pt}\includegraphics[width=0.9\linewidth]{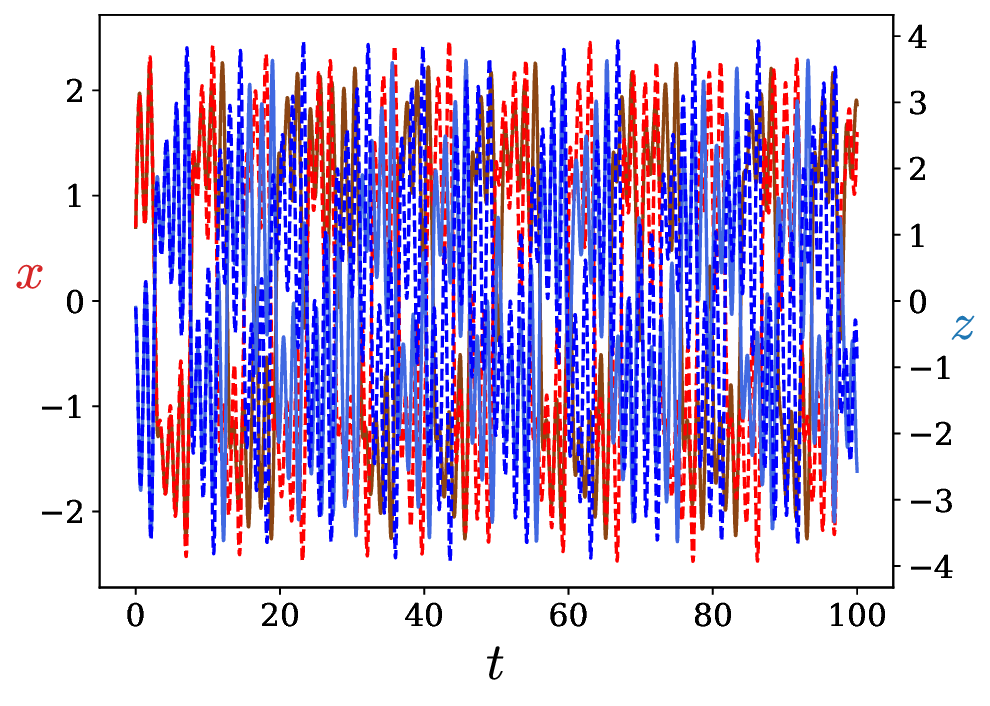}
\newline
(b) Time series
\hspace{7.5pt}\includegraphics[width=0.9\linewidth]{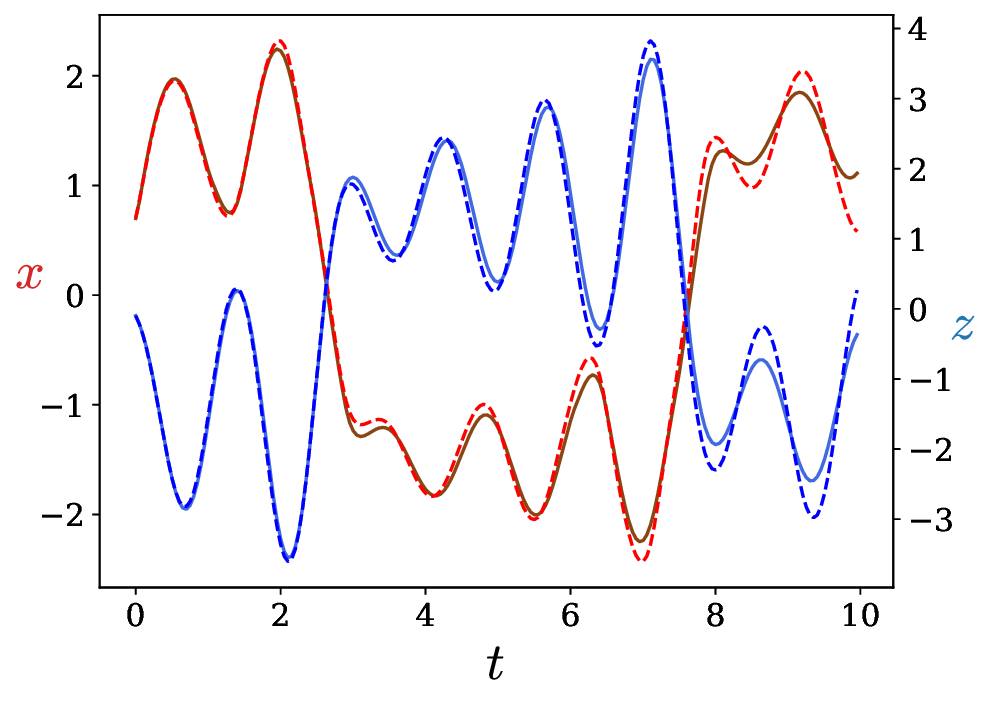}
\newline
(c) Zoomed-in view of (b)
\caption{(a) State-space comparison of the trajectory generated by Eqs.~(\ref{chua_main_x}-\ref{chua_main_z}) (dashed magenta) estimated at $\X$ (light green) with that generated by the ground truth model (Eqs.~\ref{chua_x}-\ref{chua_z}) (solid green). Without the absolute valued function in our primitives, the identified model can't give exactly the same attractor as the true model, but we can see similar dynamics and switching behavior. (b) The same trajectories plotted as a time series for $x$ and $z$ shows that while the solution of the identified model (dashed) diverges from that for the true model (solid), both show similar statistical behavior.  Moreover, the identified model tracks the true model for the initial times as shown in (c), before gradually diverging.}
\label{chua_fig}
\end{figure}

Lastly, we demonstrate the estimation of a model for a dynamical system that has functions which are not included in our operational layer. For this purpose, we consider the Chua double-scroll attractor~\cite{chua1986double}:
\begin{subequations}
\begin{align}
    \dot{x} &= 15.6 y - \frac{31.2}{7} x + 3.343 (|x + 1| - |x-1|), \label{chua_x}\\
    \dot{y} &= x - y + z,\\
    \dot{z} &= -28 y.\label{chua_z}
\end{align}
\end{subequations}
The circuit modeled by the above equations shows chaotic behavior that has been observed by an analog oscilloscope. It has one nonlinear element which is a piecewise-linear resistor. This nonlinearity is modeled by the absolute-valued function of state as $|x + 1| - |x-1|$, which is not a differentiable function of state. To test our algorithm, we again consider $1000$ data points from a single randomly initialized trajectory (equally spaced in time on the time interval [0, 100]) as shown in Fig.~\ref{chua_fig}(a). We use $L = 10$ layers and $K = 1$ stack with $322$ trainable parameters. Without including the absolute value function in the operational sublayers, we know that this model cannot be identified exactly. We find that most of the identified models consist of one or more sinusoidal terms, such as:
\begin{subequations}
\begin{align}
    \dot{x} &= 15.857~y - \frac{1}{4} x + 2.143~\sin(2~x), \label{chua_main_x}\\
    \dot{y} &= x - y + z,\\
    \dot{z} &= -28.667~y.\label{chua_main_z}
\end{align}
\end{subequations}
This model is identified in $6400$ epochs on training data $\Y$ at the input data points $\X$ with a relative-RMSE of $\approx~5.41\%$. For reference, the ground truth model gives an error of $\approx~1.83\%$ in Table~\ref{rmse_table}.  We see from Fig.~\ref{chua_fig}(c) that the selected network model gives good short-time tracking of the trajectory of the exact model, and from Fig.~\ref{chua_fig}(b) that they have similar long-time statistical behavior. SymANNTEx identified this model up to noise levels of $\sigma_1 = 0.01, \sigma_2 = 0.01$.

\section{Comparison with $L_1$ regularized Mean Squared Error as the loss}\label{comparison}

Our choice of using the Mean Absolute Error (MAE) alongside $L_{1/2}$, $L_{poly}$  and $L_{ops}$ regularization as the loss function (Eq.~\ref{loss}) to be optimized differs from the common usage of $L_1$ regularized Mean Squared Error (MSE) in the literature~\cite{mohri2018foundations}. Despite working in a Euclidean space and our own metric for performance presented above being RMSE, this choice is driven by empirical observations that the loss function of MAE with custom regularization identifies models that are far more parsimonious than those identified using $L_1$ regularized MSE as the loss function. While this could in part be attributed to the sparsity promoting $L_{1/2}$ regularization~\cite{Xu2010}, we found that its usage with MSE did not generate models as parsimonious as our choice of loss function. Some analysis and comparison between MAE and MSE as loss functions has been done~\cite{qi2020mae} in the context of conventional DNNs without regularization. We summarize preliminary findings on this in Sec.~\ref{comparison}, but detailed investigation is beyond the scope of this work. Here, we present some examples of models identified with the conventional $L_1$ regularized MSE. These examples are summarized in Table \ref{param_table}. In comparison to our custom loss function, we see that some of the identified models have far larger Relative-RMSE than the training noise itself. Deep learning literature suggests that this could arise from underfitting~\cite{goodfellow2016deep}, but we can see that the models identified using SymANNTEx have a higher number of parameters than both the ground truth models and selected network models in Table \ref{param_table}. However, there are also models with acceptable Relative-RMSE, but we can see that those models seem to have been overfitted~\cite{goodfellow2016deep} to the training data using the corresponding number of parameters, which are more than those in the ground truth models. However, we noticed that the R\"{o}ssler model was consistently identified. In our numerical investigations  using other combinations of error and regularization namely: MAE with $L_1$ regularization, MAE with $L_1, \ L_{poly}, \ L_{ops}$ regularizations, and MSE with $L_1, \ L_{poly}, \ L_{ops}$ regularizations, we see trends similar with some models being identified well while others being overfitted. Thus we avoid tabulating them in the interest of brevity and summarize only the results from $L_1$ regularized MSE in Table~\ref{param_table} as it is the most commonly used one.

We remark that the AIC score, which considers the number of non-zero parameters alongside MSE, could be seen as the so-called $L_0$ regularization. Although we do not optimize the network parameters using this in the training process, our model selection through computing AIC scores clearly takes this into account.

\begin{table*}\centering
\caption{Selected network models for demonstrating versatility of SymANNTEx. Note the similarity of the Relative-RMSEs with those in Table~\ref{rmse_table}, and the number of parameters in the estimated and ground truth models.}\label{param_table}
\begin{ruledtabular}

\begin{tabular}{cccccc}
Example & \multicolumn{2}{|c|}{Identified model} & Parameters in & \multicolumn{2}{|c}{$L_1$ regularized MSE} \\ &  \multicolumn{1}{|c}{$\approx$ Relative-RMSE} & \multicolumn{1}{c|}{Parameters} & ground truth model & \multicolumn{1}{|c}{$\approx$ Relative-RMSE} & \multicolumn{1}{c}{Parameters}\\
\hline
Takens-Bogdanov &  $2.10\%$ & 4 & 4 & $ 2.82\%$ & 10 \\
Simple Pendulum &  $4.45\%$ & 2 & 2 & $ 66.78\%$ & 3\\
R\"{o}ssler &  $2.60\%$ & 4 & 4 & $ 2.62\%$ & 4\\
Lorenz &  $2.59\%$ & 5 & 5 & $ 5.09\%$ & 54\\
FitzHugh-Nagumo &  $2.40\%$ & 10 & 8 & $ 6.00\%$ & 11\\
Chemical kinetics & $1.56\%$ & 11 & 7 & $ 0.61\%$ & 48\\
Chua Double-scroll & $5.41\%$ & 8 & 10 & $8.43\%$ & 18\\
\end{tabular}
\end{ruledtabular}
\end{table*}


\section{Discussion}
We presented a novel algorithm for identifying interpretable, closed-form models for a dynamical system from its time series data, based on simple operations on state variables and functions.  It is a generative dictionary system identification method, which only pre-specifies a set of primitive operations and functions, and creates the dictionary of possible model terms from combinations and compositions of these primitives.  Moreover, it is a symbolic regression algorithm, in that it searches a space of mathematical expressions to find the model, as opposed to conventional regression techniques which optimize parameters for a fixed dictionary.  These characteristics allow a much wider array of possible model terms than fixed dictionary system identification methods.  Our neural network architecture is novel in its utilization of the weights of the network to determine the form of the terms in the generated dictionary, the number of which remains remains small enough to handle computationally, but large enough to capture a rich set of possibilities. The depth of our network, characterized by the number of stacks, represents the level of complexity obtained by composing functions.  Each stack is made up of operational layers that span the breadth of the network, and which account for multiple occurrences of a function in the model but with different coefficients.  The AIC score was used to select between models with parameters chosen within given tolerances.  Unlike deep learning and reservoir computing approaches for forecasting a system's state~\cite{timeseriesreview, gauthier2021next}, SymANNTEx gives closed-form models for the dynamics, which could be useful for designing control algorithms.  A powerful feature of our approach is that it builds upon desired primitive functions via compositions and linear combinations rather than requiring pre-specified fixed basis functions whose linear combination describes the system's dynamics. When a requisite function is absent from our chosen primitives, say, the absolute value (which is also non-differentiable) we found that SymANNTEx gives its best estimate with possible terms.
Our algorithm gives accurate models even when the derivative data is noisy, but there can be a trade-off between the accuracy of the model and the number of terms which it contains.  It would be interesting to explore the use of recently proposed approaches for model selection with SymANNTEx when the data is more highly corrupted by noise.\cite{tran17, almo20}

We found empirically that training data with larger $\|\cdot\|_2$ norm averaged over the data-points -- a notion of ``power'' in the derivatives --  leads to more accurate models with smaller generalization error. Preliminary analysis of metrics during training indicate that data with lower power in the derivatives are prone to the occurrence of vanishing gradients~\cite{goodfellow2016deep}. Thus, as for SINDy~\cite{brun16} and EDMD~\cite{edmd}, we do not normalize the training data. This is in contrast to algorithms where normalizing training data is considered to be better practice~\cite{goodfellow2016deep}. Despite this, SymANNTEx has successfully identified the FitzHugh-Nagumo model and the model with Arrhenius terms, both of which exhibit separation of time-scales.

We observed several interesting artifacts in training.  We found that a loss function constituted by MAE alongside sublayer-specific regularization -- some of which are non-convex like the $L_{1/2}$ regularization~\cite{Xu2010} -- outperforms the conventional combination of MSE with $L_1$ regularization in estimating parsimonious models generalizable beyond seen data $\X$, even when trained for $10^5$ epochs. The generalizable aspect could also be related to the recently reported double-descent phenomenon~\cite{double_descent} where increasing the number of trainable parameters beyond the widely accepted optimal range in bias-variance trade-off seems to promote generalizability. An example is the usage of $L = 10$ layers in identifying the Takens-Bogdanov normal form, the Lorenz, and the R\"{o}ssler equations: all of these should in principle be identifiable with $L = 2$ layers but we find that two layers are consistently insufficient. We also observed faster training at the edge of stability~\cite{edge_of_stability} with initial learning rates higher than the default in Adam~\cite{kingma2015adam}. For example, training using default learning rate constant of $0.001$ even for over $20,000$ epochs yields residual terms in the identified Lorenz and FitzHugh-Nagumo equations that vanish in as few as $800$ epochs with learning rate constant of $0.01$.

\section*{Acknowledgement}
This work was supported by National Science Foundation Grant No. NSF-2016004.

\nocite{*}
\bibliography{aipsamp}

\end{document}